\documentclass[11pt,twoside]{amsart}
\usepackage[dvipsnames]{xcolor}
\usepackage{hyperref}
\hypersetup{colorlinks=true, linkcolor=Blue, citecolor=Maroon}
\usepackage{amsthm, amsmath, amscd, amssymb,centernot,txfonts}
\usepackage{tikz-cd}
\usepackage{lipsum}
\usepackage{cancel} 
\usepackage{multicol}
\usepackage{amsfonts}
\usepackage{mathrsfs}
\usepackage[all]{xy}
\usepackage[left=2.5cm,top=2.5cm,bottom=3cm,right=2.5cm]{geometry}
\usepackage{dirtytalk}
\usepackage{mathtools}
\usetikzlibrary{shapes,arrows}
\usepackage{etoolbox}
\usepackage{dynkin-diagrams}
\usepackage{enumitem}
\usepackage{chngpage}
\usepackage{adjustbox}
\usepackage[normalem]{ulem}
\usepackage[utf8]{inputenc}
\usepackage{fourier} 
\usepackage{array}
\usepackage{makecell}
\usepackage{cleveref}
\usepackage{comment}
\usepackage{cancel}
\usepackage{caption}
\newcommand{\longsquiggly}{\xymatrix{{}\ar@{~>}[r]&{}}}
\usepackage{graphicx}

\theoremstyle{plain}
\numberwithin{equation}{section}
\newtheorem{theorem}{Theorem}[section]
\newtheorem{proposition}[theorem]{Proposition}
\newtheorem{lemma}[theorem]{Lemma}
\newtheorem{corollary}[theorem]{Corollary}

\newtheorem{metaprinciple}[theorem]{Metaprinciple}

\theoremstyle{definition}
\newtheorem{remark}[theorem]{Remark}

\newtheorem{set-up}[theorem]{Set-up}
\newtheorem{definition}[theorem]{Definition}

\newcommand*{\QEDA}{\hfill\ensuremath{\blacksquare}}

\tikzstyle{decision} = [diamond, draw, , 
    text width=4.5em, text badly centered, node distance=3cm, inner sep=0pt]
\tikzstyle{block} = [rectangle, draw, , 
    text width=10em, text centered, rounded corners, minimum height=2em]
\tikzstyle{block1} = [rectangle, draw, , 
    text width=5em, text centered, rounded corners, minimum height=2em]    
\tikzstyle{line} = [draw, -latex']
\tikzstyle{cloud} = [draw, ellipse,, node distance=3cm,
    minimum height=2em]
\begin{document}

\title[Extendability of projective varieties via degeneration to ribbons]{Extendability of projective varieties via degeneration to ribbons with applications to Calabi-Yau threefolds}

\author[P. Bangere]{Purnaprajna Bangere}
\address{Department of Mathematics, University of Kansas, Lawrence, USA}
\email{purna@ku.edu}


\author[J. Mukherjee]{Jayan Mukherjee}
\address{Department of Mathematics, Oklahoma State University, Stillwater, USA}
\email{jayan.mukherjee@okstate.edu}


\maketitle
\begin{abstract}
In this article we study the extendability of a smooth projective variety by degenerating it to a ribbon. We apply the techniques to study extendability of Calabi-Yau threefolds $X_t$ that are general deformations of Calabi-Yau double covers of Fano threefolds of Picard rank $1$. The Calabi-Yau threefolds $X_t \xhookrightarrow[]{|lA_t|} \mathbb{P}^{N_l}$, with $l \geq j$, where $A_t$ is the generator of Pic$(X_t)$ and $j$ is the least positive integer such that $lA_t$ is very ample, are general elements of a unique irreducible component $\mathcal{H}_l^Y$ of the Hilbert scheme which contains embedded Calabi-Yau ribbons on $Y$ as a special locus. For $l = j$, using the classification of Mukai varieties, we show that the general Calabi-Yau threefold parameterized by $\mathcal{H}_j^Y$ is as many times smoothly extendable as $Y$ itself. On the other hand, we find for each deformation type $Y$, an effective integer $l_Y$ such that for $l \geq l_Y$, the general Calabi-Yau threefold parameterized by $\mathcal{H}_l^Y$ is not extendable. These results provide a contrast and a parallel with the lower dimensional analogues; namely, $K3$ surfaces and canonical curves (\cite{CLM93}, \cite{CLM93A}, \cite{CLM98}) which stems from the following result we prove: for $l \geq l_Y$, the general hyperplane sections of elements of $\mathcal{H}_l^Y$ fill out an entire irreducible component $\mathcal{S}_l^Y$ of the Hilbert scheme of canonical surfaces which are precisely $1-$ extendable with $\mathcal{H}^Y_l$ being the unique component dominating $\mathcal{S}_l^Y$. The contrast lies in the fact that for polarized $K3$ surfaces of large degree, the canonical curve sections do not fill out an entire component while the parallel is in the fact that the canonical curve sections are exactly one-extendable.

\end{abstract}

\section{Introduction}

\begin{definition}\label{extendability}
Let $X \subset \mathbb{P}^N$ be an irreducible nondegenerate variety of codimension at least $1$. Let $k \geq 1$ be an integer. We say that $X$ is $k-$ extendable if there exists a variety $W \subset \mathbb{P}^{N+k}$
different from a cone, with dim $W =$ dim $X + k$ and having $X$ as a section by a $N-$ dimensional linear space such that $W$ is smooth along $X = W \cap \mathbb{P}^N$. We say that $X$ is precisely $k-$ extendable if it is $k-$ extendable but not $(k + 1)-$ extendable. The variety $W$ is called a $k$-extension of $X$. We say that $X$ is extendable if it is $1-$ extendable.  
\end{definition}

The extendability of a projective variety is a natural and fundamental question in projective geometry. This classical question has been the topic of intense research for decades and has revealed deep connections between the geometry of the embedding, Gaussian map of curve sections, deformations of cones over the hyperplane sections etc. We refer to \cite{Lop23} for an excellent recent survey of this topic. 
The purpose of this article is to introduce techniques to study extendability of smooth projective varieties via degeneration to ribbons. In particular we show two different themes, one to show extendability and the other to show non-extendability of the general member of a one parameter family that projectively degenerates to a ribbon in the central fiber. \textcolor{black}{We apply them to study extendability of Calabi-Yau threefolds, where not much is known.} \textcolor{black}{It is interesting to note that ribbons have appeared in the context of extendability of canonical curves and $K3$ surfaces in \cite{CDS20}, \cite{V92}, but in a way that is completely different from this article.} \textcolor{black}{We will now recall the definition of a ribbon}

\begin{definition}
   A ribbon $\widetilde{Y}$ on a reduced connected scheme $Y$ is an everywhere non-reduced scheme with $\widetilde{Y}_{\textrm{red}} = Y$ such that
   \begin{itemize}
       \item[(1)] The ideal sheaf $\mathcal{I}_{Y/\widetilde{Y}}$ of $Y$ inside $\widetilde{Y}$ satisfies $\mathcal{I}_{Y/\widetilde{Y}}^2 = (0)$ and 
       \item[(2)] $\mathcal{I}_{Y/\widetilde{Y}}$ as an $\mathcal{O}_Y$ module is a locally free sheaf $L$ of rank one also called the conormal bundle of the ribbon.
   \end{itemize}

\subsection{Extendability of the general member of a component of the Hilbert scheme  containing ribbons} We now explain how one can use ribbon degenerations to derive information on the extendability of a general member of a Hilbert component. Suppose one has an irreducible component of the Hilbert scheme $\mathcal{H}$ of $\mathbb{P}^N$ which parameterizes smooth varieties in general and ribbons $\widetilde{Y}$ on $Y$ with conormal bundle $L$ as a special locus which represents smooth points of $\mathcal{H}$. Now suppose that we have the knowledge of extendability of the induced embedding of the reduced part $Y \subset \mathbb{P}^M$ ($M \leq N$), i.e, $W \subset \mathbb{P}^{M+k}$ is a smooth $k-$ step extension of $Y \subset \mathbb{P}^M$. Then we can derive knowledge of the extendability of the general elements of $\mathcal{H}$. \textcolor{black}{We start by showing that the smooth members in a certain special locus in $\mathcal{H}$ are $k-$ extendable} as follows :

\begin{enumerate}
    \item Extend the ribbon $\widetilde{Y} \subset \mathbb{P}^{N}$ to a ribbon $\widetilde{W} \subset \mathbb{P}^{N+k}$. This eventually boils down to showing existence of ribbons supported on $W$ 
    in $\mathbb{P}^{N+k}$ extending the embedding $W \subset \mathbb{P}^{M+k} \subset \mathbb{P}^{N+k}$ which is parameterized by nowhere vanishing sections of $H^0(N_{W/\mathbb{P}^{N+k}} \otimes \mathcal{L})$ where $\mathcal{L}$ is a line bundle on $W$ whose restriction to $Y$ is $L$.
    
    \item Show that the extended ribbon $\widetilde{W} \subset \mathbb{P}^{N+k}$ is smoothable.
    Thus we have shown the $k-$ step extendability of some smooth \textcolor{black}{special} members of $\mathcal{H}$ containing the ribbons in their closure.

    \item \textcolor{black}{Finally we conclude by the showing that} the corresponding \textcolor{black}{non-empty} flag Hilbert scheme \textcolor{black}{dominates} $\mathcal{H}$.

\end{enumerate}

We remark that the abstract variety $Y$ or the embedding of $Y \subset \mathbb{P}^M$ is usually a \textcolor{black}{cohomologically} simpler variety \textcolor{black}{with a} much better understood embedding than $X \subset \mathbb{P}^N$ where $X$ is a general member of $\mathcal{H}$. 
\textcolor{black}{So}
it is reasonable to expect that we will have better knowledge about its extendability.  
\textcolor{black}{In the context of this article, we show the extendability of general smoothings of Calabi-Yau ribbons on anticanonically embedded Fano-threefolds of Picard rank one. Guided by the metaprinciple \ref{extendability}, we accomplish this by showing the existence of very specific Fano ribbon structures on Mukai varieties, which are extensions of the three dimensional Calabi Yau ribbons, and then by smoothing them.}

\subsection{Non-extendability of the general member of a component of the Hilbert scheme  containing ribbons} To show non-extendability, we apply Zak-L\'vovsky's criterion on the ribbon $\widetilde{Y} \subset \mathbb{P}^N$. We recall the definition of the invariant $\alpha(X)$ and two results by Zak and L\'vovsky (\cite{Lvo92}, \cite{Zak91}, see also \cite{Bad94}, \cite[Theorem $1.3$]{Lop23}).

\begin{definition}

For $X \subset \mathbb{P}^{N}$ a smooth irreducible nondegenerate variety of codimension at least $1$ with normal bundle $N_{X/\mathbb{P}^N}$, we set
$\alpha(X) = h^0(N_{X/\mathbb{P}^N}(-1))-N-1$

\end{definition}

\begin{theorem}(\cite{Lvo92}, \cite{Zak91})\label{ZL}
Let $X \subset \mathbb{P}^{N}$ be a smooth irreducible nondegenerate variety of codimension at least $1$ and suppose $X$ is not a quadric. If $\alpha(X) \leq 0$, then X is not extendable. Further, given an integer $k \geq 2$, suppose that either:
\begin{itemize}
    \item[(1)] $\alpha(X) < N$ or
    \item[(2)] $H^0(N_{X/\mathbb{P}^N}(-2)) = 0$,
\end{itemize}
If $\alpha(X) \leq k-1$, then X is not $k-$ extendable.  

\end{theorem}

\end{definition}
\textcolor{black}{By results in \cite{BE95}, \cite{Fon93}, \cite{Gon06}, \cite{GGP10}, \cite{GGP13}, it is known that under sufficiently general conditions, as an embedding degenerates to a double cover, the embedded models degenerate to a ribbon over the base of the double cover. \textcolor{black}{Since the ribbon is a local complete intersection,} by a series of short exact sequences, it is shown that computing the cohomology of the twists of the normal bundle    $N_{\widetilde{Y}/\mathbb{P}^N}$ of a embedded ribbon $\widetilde{Y}$ on $Y\subset \mathbb{P}^N$ eventually boils down to computing the cohomology of the twists of the normal bundle $N_{Y/\mathbb{P}^N}$ (see equations \ref{1}, \ref{2}, \ref{3}) which are easier to handle than the normal bundle of a general member $X \subset \mathbb{P}^N$ in the family, about which we know very little.}
 
\textcolor{black}{In this paper we study the extendability of Calabi-Yau threefolds $X \subset \mathbb{P}^N$ by degenerating \textcolor{black}{them} to Calabi-Yau ribbons on Fano-threefolds $Y$ of Picard rank one, embedded by the complete linear series of $-lH$ where $H$ is the generator of Pic$(Y)$. For $l$ equal to the index of $Y$, we show that the Calabi-Yau threefolds are at least as many times smoothly extendable as the Fano threefold $Y$. For higher values of $l$, we show non-extendability of $X$.} Fano threefolds of Picard rank $1$ are classified into $17$ different deformation types by the classification of Iskovskih-Mori-Mukai (\cite{I77}, \cite{I78}, \cite{MM82}). We use the explicit description of each of the families which are either complete intersections in weighted projective spaces or regular sections of homogeneous vector bundles on Grassmannians (see \cite{DFT22}, \cite{Bel24}, \url{https://www.fanography.info/.}). The computations boil down to Borel-Bott-Weil theorems. We now state our results. First we state a Proposition describing the Calabi-Yau threefolds we work with and their Hilbert schemes.

\begin{proposition}\label{Intro Prop 1}(see Proposition \ref{existence and smoothing of ribbons}, Remark \ref{values of j})
 Let $Y$ be a smooth, projective Fano threefold of Picard rank $1$. Let $i(Y)$ denote the index of $Y$. Let $\pi: X \to Y$ denote a Calabi-Yau double cover branched along a smooth divisor $B \in |-2K_Y|$. Let $H$ denote the generator of Pic $(Y)$. Let $j\geq i(Y)$ be the least positive integer such that $jH$ is very ample.Then
\begin{enumerate}
     \item Pic $(X) = \mathbb{Z}$. If $A$ denotes the generator of Pic $(X)$, then for any polarized deformation $(X_t, A_t)$ of $(X,A)$ along a one parameter family $T$ the line bundle $A_t$ generates Pic $(X_t)$. 
     \item \textcolor{black}{For $l \geq j$, $h^0(lA_t) = h^0(lH)+h^0((l-i(Y))H)$. We set $N_l+1 = h^0(lH)+h^0((l-i(Y)H)$.}
     \item \textcolor{black}{Let $l \geq j$. If $Y$ belongs to family $1.11$, further assume $l \geq 4$. Then for a general polarized deformation $(X_t, A_t)$ of $(X,A)$ along a one parameter family $T$, $lA_t$ is very ample. Further there exists a unique irreducible component $\mathcal{H}^Y_l$ of the Hilbert scheme of $\mathbb{P}^{N_l}$ with Hilbert polynomial $p(z) = h^0(lzH)+h^0((lz-i(Y))H)$ that parameterizes linearly normal Calabi-Yau varieties $X_t \subset \mathbb{P}^{N_l}$ such that $(X_t, \mathcal{O}_{X_t}(1))$ can be deformed to $(X, lA)$ along an irreducible one parameter family.}
     \item The component $\mathcal{H}^Y_l$ contains linearly normal embedded Calabi-Yau ribbons $\widetilde{Y} \subset \mathbb{P}^{N_l}$ on $Y$ such that $\mathcal{O}_{\widetilde{Y}}(1)|_Y = lH$. \textcolor{black}{The embedded ribbons form an irreducible family, and they represent smooth points} of $\mathcal{H}^Y_l$. \textcolor{black}{Consequently, any linearly normal embedded Calabi-Yau ribbon $\widetilde{Y} \subset \mathbb{P}^{N_l}$ on $Y$ such that $\mathcal{O}_{\widetilde{Y}}(1)|_Y = lH$ is contained in $\mathcal{H}^Y_l$}

     \item \textcolor{black}{If $Y$ is not of type $1.1$, $lA_t$ is not very ample for $l < j$.}  
     
\end{enumerate}

\end{proposition} 

With the notation as above, we now state our main result on the extendability of general Calabi-Yau threefolds in $\mathcal{H}^Y_l$ : 

\begin{theorem}(see Theorem \ref{l=j}, Theorem \ref{Picard rank one}) \label{intro Thm 1}
    \begin{itemize}
        \item[(1)] The extendability for a general smooth Calabi-Yau threefold of the Hilbert component $\mathcal{H}^Y_l$ ($l \geq j$) (in the notation of Proposition \ref{existence and smoothing of ribbons}) is listed in the fourth columns of tables \hyperref[t01]{1}, \hyperref[t02]{2}. \textcolor{black}{The tables show that out of the $17$ different deformation types of $Y$,} 
        \begin{enumerate}
            \item[(a)] \textcolor{black}{Excepting for families $1.1$ and $1.11$, when $l = j$, a smooth Calabi-Yau threefold $X_t \subset \mathbb{P}^{N_l}$ for a general $t$ in $\mathcal{H}^Y_j$ is at least as much extendable as the anticanonically embedded Fano $Y$.}

            \item[(b)] \textcolor{black}{For each family, we give an effective value $l_Y$, such that for $l \geq l_Y$, a smooth Calabi-Yau threefold $X_t \subset \mathbb{P}^{N_l}$ for a general $t$ in $\mathcal{H}^Y_l$ has $\alpha(X_t) = 0$ and is hence not extendable. Consequently, the cone over $X_t$ in $\mathbb{P}^{N_l+1}$ is not smoothable in $\mathbb{P}^{N_l+1}$.}
        \end{enumerate}
        \item[(2)] For each deformation type $Y$ and $l \geq l_Y$, the general canonical surface sections $S_t \subset \mathbb{P}^{N_l-1}$ of the \textcolor{black}{general} Calabi-Yau threefolds $X_t \subset \mathbb{P}^{N_l}$ of $\mathcal{H}^Y_l$ have $\alpha(S_t) = 1$, are exactly $1-$ extendable and fill out a unique irreducible component of the Hilbert scheme $\mathcal{S}^Y_l$ with a dominant map $\mathcal{H}^Y_l \to \mathcal{S}^Y_l$. 
        \item[(3)] For $l \geq l_Y$, $\mathcal{H}^Y_l$ is the only irreducible component of the Hilbert scheme parameterizing linearly normal Calabi-Yau threefolds that dominate the component $\mathcal{S}^Y_l$ by taking hyperplane sections. For a general canonical surface $S_t \subset \mathbb{P}^{N_l-1} \in \mathcal{S}^Y_l$, the Calabi-Yau threefolds in $\mathbb{P}^{N_l}$ containing $S_t$ as a hyperplane section forms an irreducible family. 
    \end{itemize}
\end{theorem}

It is interesting to compare the results of this article with the results on extendability of $K3$ surfaces and canonical curves in \cite{CLM93}, \cite{CLM93A} and \cite{CLM98} (see also \cite{CD22}). Let $\mathcal{H}_{r,g}$ denote the Hilbert scheme of $K3$ surfaces embedded by the $r-$ th veronese map of a primitive very ample line bundle of sectional genus $g$. By studying Gaussian maps of canonical curve sections of such $K3$ surfaces, it was proven that for $(r = 1, g \geq 13), (r = 2, g \geq 7), (r = 3, g \geq 5), (r = 4, g \geq 4), (r \geq 5, g \geq 3)$ (see \cite[Table $2.14$]{CLM98} for the precise and slightly stronger results), the general curve section of the general $K3$ surfaces in $\mathcal{H}_{r,g}$ has corank one Gaussian map which in particular means that the general $K3$ surfaces are not extendable and the general curve sections are exactly one extendable. The results Theorem \ref{intro Thm 1}, $(1)$, $(3)$ are analogues of the above results. However, there is a \textcolor{black}{contrast for Calabi-Yau threefolds}; the hyperplane sections of $\mathcal{H}^Y_l$ for $l \geq l_Y$ fill out an entire irreducible component of the Hilbert scheme of canonical surfaces, while curves on a $K3$ surface are non-general curves of the corresponding component of the Hilbert scheme. We plot the invariants of the surface sections in subsection \ref{invariants of hyperplane sections}. 

\noindent\textbf{Acknowledgements.} We thank Professor Ciro Ciliberto for his encouragement. \textcolor{black}{We thank the anonymous referee for carefully reading our paper. The referee's comments and suggestions, helped us revise some computations and improve the clarity and exposition of the paper.}

\section{Main result on extendability of projective varieties }\label{main results}

The following general metaprinciple forms the basis to show extendability of general embedded smoothings of ribbons using the knowledge of extendability of the reduced part of the ribbon. \textcolor{black}{By a one-parameter family we mean a flat projective morphism $\mathcal{X} \to T$ where $T$ is a smooth irreducible curve.} 

\begin{metaprinciple}\label{extendability}

\begin{enumerate}
    \item Let $Y \subset \mathbb{P}^M$ be an embedding of a smooth projective variety. Let $W \subset \mathbb{P}^{M+k}$ be a smooth extension of $Y \subset \mathbb{P}^M$. Suppose that there exists ribbons $\widetilde{W} \subset \mathbb{P}^{N+k}$ with conormal bundle $\widetilde{L}$ such that $\widetilde{L}|_Y = L$, extending the embedding $W \subset \mathbb{P}^{M+k} \subset \mathbb{P}^{N+k}$, where the second embedding is an embedding of linear spaces. Assume further that $\widetilde{W}$ is smoothable inside $\mathbb{P}^{N+k}$ along a one-parameter family $T$. Then by taking a codimension $k$ linear section of the smoothing family, we see that there exist embedded ribbons $\widetilde{Y} = \widetilde{W} \bigcap \mathbb{P}^N \subset \mathbb{P}^N$ extending $Y \subset \mathbb{P}^M$ with conormal bundle $L$ which are smoothable along a one-parameter family $T$. Further for a general $t \in T$, the smoothing fiber $V_t \subset \mathbb{P}^{N+k}$ of $\widetilde{W} \subset \mathbb{P}^{N+k}$ is a $k-$ step extension of the smoothing fiber $X_t \subset \mathbb{P}^N$ of $\widetilde{Y} \subset \mathbb{P}^N$.

    \item Consequently, in a situation where we are interested in the extendability of a general element of an irreducible component $\mathcal{H}$ of the Hilbert scheme in $\mathbb{P}^N$ which contains ribbons $\widetilde{Y}$ on $Y$ with conormal bundle $L$, such that the ribbons represent smooth points of $\mathcal{H}$, we can use part one to deduce the extendability of at least a special locus parameterizing smooth varieties in $\mathcal{H}$.

    \item Ribbons $\widetilde{Y} \subset \mathbb{P}^N$ on $Y$ with conormal bundle $L$ extending a given (possibly degenerate) embedding $Y \subset \mathbb{P}^N$ are parameterized by the nowhere vanishing sections of $H^0(N_{Y/\mathbb{P}^N} \otimes L)$ (see \cite[Proposition $2.1$]{Gon06}).
\end{enumerate}

\end{metaprinciple}

In the following theorem, we show how to compute an upper bound to $\alpha(\widetilde{Y})$ of an embedded ribbon $\widetilde{Y} \subset \mathbb{P}^N$ in terms of the cohomology of the twists of the normal bundle and eventually that of the tangent bundle of its reduced structure $Y$. \par

\smallskip

\noindent \textcolor{black}{We introduce the following notation that will be used in the subsequent sections.}

\smallskip

\color{black}
\noindent \textbf{Notation:} 
\begin{enumerate} 
    \item[(a)]  Let $\widetilde{V}$ is a vector bundle on the ribbon $\widetilde{Y}$ and $L$ is a line bundle on $Y$. Then $L$ has the structure of an $\mathcal{O}_{\widetilde{Y}}$ module due to the homomorphism $\mathcal{O}_{\widetilde{Y}} \to \mathcal{O}_Y$ that defines the ribbon $\widetilde{Y}$.  We define $\widetilde{V}(L):=\widetilde{V}\otimes_{\mathcal{O}_{\widetilde{Y}}} L$.

    \item[(b)] By $I_{\widetilde{Y}}$ and $I_{Y}$, we mean the ideal sheaf of $\widetilde{Y}$ and $Y$ in $\mathbb{P}^N$ respectively.
\end{enumerate}
  
\color{black}

\begin{theorem}\label{ZL for ribbons}
 Let $\widetilde{Y}$ be a ribbon over a smooth projective variety $Y$ of dimension $d > 1$ with conormal bundle $L$. Let $\widetilde{Y} \hookrightarrow \mathbb{P}^N$ be an embedding of $\widetilde{Y}$ induced by the complete linear series of a very ample line bundle $\widetilde{H}$. Let $H = \widetilde{H}|_Y$ and assume $H^1(-H+L) = 0$. Let
 $$\beta = h^1(T_Y(-H+L)) + h^1(T_Y(-H)) - h^0(T_Y(-H)) +  h^1(-H-L) + h^0(-H-2L)$$ 
Then 
$$\alpha(\widetilde{Y}) = h^0(N_{\widetilde{Y}/\mathbb{P}^N}(-\widetilde{H}))-N-1 \leq  \beta$$ 
 
\begin{enumerate}
    \item If $\beta < N$ and $\widetilde{Y}$ is smoothable inside $\mathbb{P}^N$ into smooth fibers $X_t \subset \mathbb{P}^N$, then a general $X_t$ is not $\beta+1$ extendable.

   \item If further, $\beta = 0$, then $\alpha(\widetilde{Y}) \leq 0$ (or $h^0(N_{\widetilde{Y}/\mathbb{P}^N}(-\widetilde{H})) \leq N+1$). Consequently, if $\widetilde{Y}$ is smoothable inside $\mathbb{P}^N$ into smooth fibers $X_t \subset \mathbb{P}^N$, then, $\alpha(X_t) = 0$ and a general $X_t$ is not extendable. In this case, $H^0(N_{X_t/\mathbb{P}^N}(-2)) = 0$ and if $H^1(\mathcal{O}_{X_t}) = 0$ (which is true if $H^1(\mathcal{O}_Y) = H^1(L) = 0$), we have cork($\Phi_{\mathcal{O}_{C_t}(1), K_{C_t}}) \geq d-1$ for a general curve section $C_t$ of $X_t$.

\end{enumerate}

\end{theorem}

\noindent \textit{Proof.} We have the following exact sequences (see \cite[Lemma $4.2$]{GGP08})
\begin{equation}\label{1}
    0 \to N_{\widetilde{Y}/\mathbb{P}^N}(L) \to N_{\widetilde{Y}/\mathbb{P}^N} \to N_{\widetilde{Y}/\mathbb{P}^N} \otimes \mathcal{O}_Y \to 0
\end{equation}
\begin{equation}\label{2}
    0 \to \mathcal{H}om(I_{\widetilde{Y}}/I_{Y}^2, \mathcal{O}_Y) \to N_{\widetilde{Y}/\mathbb{P}^N} \otimes \mathcal{O}_Y \to \mathcal{O}_Y(-2L) \to 0 
\end{equation}
\begin{equation}\label{3}
    0 \to \mathcal{O}_Y(-L) \to N_{Y/\mathbb{P}^N} \to \mathcal{H}om(I_{\widetilde{Y}}/I_{Y}^2, \mathcal{O}_Y) \to 0
\end{equation}
Tensoring equation \ref{1} by $\mathcal{O}_{\widetilde{Y}}(-\widetilde{H})$, we have that 

\begin{equation}\label{4}
 h^0(N_{\widetilde{Y}/\mathbb{P}^N}(-\widetilde{H})) \leq h^0(N_{\widetilde{Y}/\mathbb{P}^N}(-H+L)) + h^0(N_{\widetilde{Y}/\mathbb{P}^N}(-H))
\end{equation}
We compute $h^0(N_{\widetilde{Y}/\mathbb{P}^N}(-H+L))$. Tensoring equation \ref{2} by $\mathcal{O}_Y(-H+L)$ we have 

\begin{equation*}
    0 \to \mathcal{H}om(I_{\widetilde{Y}}/I_{Y}^2, \mathcal{O}_Y)(-H+L) \to N_{\widetilde{Y}/\mathbb{P}^N}(-H+L) \to \mathcal{O}_Y(-H-L) \to 0 
\end{equation*}
Hence 

\begin{equation}\label{6}
    h^0(N_{\widetilde{Y}/\mathbb{P}^N}(-H+L)) \leq h^0(\mathcal{H}om(I_{\widetilde{Y}}/I_{Y}^2, \mathcal{O}_Y)(-H+L)) + h^0(\mathcal{O}_Y(-H-L))
\end{equation}
Tensoring equation \ref{3}, by $\mathcal{O}_Y(-H+L)$ we have 

\begin{equation*}
    0 \to \mathcal{O}_Y(-H) \to N_{Y/\mathbb{P}^N}(-H+L) \to \mathcal{H}om(I_{\widetilde{Y}}/I_{Y}^2, \mathcal{O}_Y)(-H+L) \to 0
\end{equation*}
Since $H^1(-H) = 0$, by Kodaira vanishing theorem, and $h^0(-H) = 0$, we have that 

\begin{equation}\label{7}
   h^0(\mathcal{H}om(I_{\widetilde{Y}}/I_{Y}^2, \mathcal{O}_Y)(-H+L)) = h^0(N_{Y/\mathbb{P}^N}(-H+L))  
\end{equation}
We now compute $h^0(N_{\widetilde{Y}/\mathbb{P}^N}(-H))$. Tensoring equation \ref{2} by $\mathcal{O}_Y(-H)$ we have 

\begin{equation*}
    0 \to \mathcal{H}om(I_{\widetilde{Y}}/I_{Y}^2, \mathcal{O}_Y)(-H) \to N_{\widetilde{Y}/\mathbb{P}^N}(-H) \to \mathcal{O}_Y(-H-2L) \to 0 
\end{equation*}
Hence 

\begin{equation}\label{7}
    h^0(N_{\widetilde{Y}/\mathbb{P}^N}(-H)) \leq h^0(\mathcal{H}om(I_{\widetilde{Y}}/I_{Y}^2, \mathcal{O}_Y)(-H)) + h^0(\mathcal{O}_Y(-H-2L))
\end{equation}
Tensoring equation \ref{3}, by $\mathcal{O}_Y(-H)$ we have 

\begin{equation*}
    0 \to \mathcal{O}_Y(-H-L) \to N_{Y/\mathbb{P}^N}(-H) \to \mathcal{H}om(I_{\widetilde{Y}}/I_{Y}^2, \mathcal{O}_Y)(-H) \to 0
\end{equation*}
Then we have that

\begin{equation}\label{7}
   h^0(\mathcal{H}om(I_{\widetilde{Y}}/I_{Y}^2, \mathcal{O}_Y)(-H)) \leq h^0(N_{Y/\mathbb{P}^N}(-H)) +  h^1(-H-L) - h^0(-H-L) 
\end{equation}
Hence we have that 

\begin{equation*}
 h^0(N_{\widetilde{Y}/\mathbb{P}^N}(-\widetilde{H})) \leq h^0(-H-L) + h^0(N_{Y/\mathbb{P}^N}(-H+L)) + h^0(N_{Y/\mathbb{P}^N}(-H)) +  h^1(-H-L) - h^0(-H-L) + h^0(\mathcal{O}_Y(-H-2L))
\end{equation*}

\begin{equation}\label{8}
h^0(N_{\widetilde{Y}/\mathbb{P}^N}(-\widetilde{H})) \leq  h^0(N_{Y/\mathbb{P}^N}(-H+L)) + h^0(N_{Y/\mathbb{P}^N}(-H)) +  h^1(-H-L) + h^0(\mathcal{O}_Y(-H-2L))
\end{equation}
We now compute $h^0(N_{Y/\mathbb{P}^N}(-H))$. We have

\begin{equation}\label{9}
   0 \to T_Y \to T_{\mathbb{P}^N}|_Y \to N_{Y/\mathbb{P}^N} \to 0 
\end{equation}
Tensoring equation \ref{9} by $\mathcal{O}_Y(-H)$, we have 

\begin{equation*}
    0 \to H^0(T_Y(-H)) \to H^0(T_{\mathbb{P}^N}|_Y(-H)) \to H^0(N_{Y/\mathbb{P}^N}(-H)) \to H^1(T_Y(-H))
\end{equation*}
Hence 

\begin{equation*}\label{10}
    h^0(N_{Y/\mathbb{P}^N}(-H)) \leq h^0(T_{\mathbb{P}^N}|_Y(-H)) + h^1(T_Y(-H)) - h^0(T_Y(-H))
\end{equation*}
Consider the Euler exact sequence restricted to $Y$

\begin{equation}\label{11}
   0 \to \mathcal{O}_Y \to \mathcal{O}_Y(H)^{N+1} \to T_{\mathbb{P}^N}|_Y \to 0 
\end{equation}
Tensoring equation \ref{11} by $\mathcal{O}_Y(-H)$, we have 

\begin{equation*}
    0 \to H^0(\mathcal{O}_Y)^{N+1} \to H^0(T_{\mathbb{P}^N}|_Y(-H)) \to 0 
\end{equation*}
so we have

\begin{equation}\label{12}
    h^0(T_{\mathbb{P}^N}|_Y(-H)) = N+1
\end{equation}
This implies 
\begin{equation}\label{13}
  h^0(N_{Y/\mathbb{P}^N}(-H)) \leq N+1 + h^1(T_Y(-H)) - h^0(T_Y(-H))  
\end{equation}
We now compute $h^0(N_{Y/\mathbb{P}^N}(-H+L))$. 
Tensoring equation \ref{11} by $\mathcal{O}_Y(-H+L)$, we have 

\begin{equation*}
     H^0(\mathcal{O}_Y(L))^{N+1} \to H^0(T_{\mathbb{P}^N}|_Y(-H+L)) \to H^1(\mathcal{O}_Y(-H+L)) 
\end{equation*}
Since $-L$ has a section, we have that $H^0(\mathcal{O}_Y(L)) = 0$ and by assumption $H^1(\mathcal{O}_Y(-H+L)) = 0$ and hence $H^0(T_{\mathbb{P}^N}|_Y(-H+L)) = 0$. 

Tensoring equation \ref{9} by $\mathcal{O}_Y(-H+L)$, we have 

\begin{equation*}
    H^0(T_{\mathbb{P}^N}|_Y(-H+L)) \to H^0(N_{Y/\mathbb{P}^N}(-H+L)) \to H^1(T_Y(-H+L))
\end{equation*}
Since $H^0(T_{\mathbb{P}^N}|_Y(-H+L)) = 0$, 

\begin{equation}\label{14}
    h^0(N_{Y/\mathbb{P}^N}(-H+L)) \leq h^1(T_Y(-H+L))
\end{equation}
Plugging everything back into equation \ref{8}, we have

\begin{equation}\label{15}
h^0(N_{\widetilde{Y}/\mathbb{P}^N}(-\widetilde{H})) \leq  N+1 + h^1(T_Y(-H+L)) + h^1(T_Y(-H)) - h^0(T_Y(-H)) +  h^1(-H-L) + h^0(-H-2L)
\end{equation} 

\QEDA

In the following lemma, we find out the vanishing condition to ensure that for a hyperplane section $Z \subset \mathbb{P}^{N-1}$, of $X \subset \mathbb{P}^N$, we have $\alpha(Z) = \alpha(X)+1$.

\begin{lemma}\label{ZL for surface sections}
Let $X \subset \mathbb{P}^{N}$ be a smooth irreducible nondegenerate variety of dimension $d \geq 2$. Let $Z \subset \mathbb{P}^{N-1}$ be a smooth hyperplane section. Let $\mathcal{O}_X(1)$ and $\mathcal{O}_Z(1)$ denote the hyperplane of $X$ and $Z$ repectively. Then 
\begin{enumerate}
    \item If $H^0(N_{X/\mathbb{P}^N}(-2)) = 0$, $\alpha(Z) \geq \alpha(X)+1$. Consequently 
    \item If $H^0(N_{X/\mathbb{P}^N}(-2)) = 0$ and $H^1(N_{X/\mathbb{P}^N}(-2)) \to H^1(N_{X/\mathbb{P}^N}(-1))$ is injective, then $\alpha(Z) = \alpha(X)+1$ 
\end{enumerate}

\end{lemma}

\noindent\textit{Proof.} First note that $N_{X/\mathbb{P}^N}|_Z = N_{Z/\mathbb{P}^{N-1}}$. Then for any $i$ we have the exact sequence
$$0 \to N_{X/\mathbb{P}^N}(-i-1) \to N_{X/\mathbb{P}^N}(-i) \to N_{Z/\mathbb{P}^{N-1}}(-i) \to 0 $$

Then for $i = 1$, we have $\alpha(X) = h^0(N_{X/\mathbb{P}^N}(-1))-N-1 \leq h^0(N_{Z/\mathbb{P}^{N-1}})(-1)-N-1 = \alpha(Z)-1$. Further from the same exact sequence, we have the equality under the added hypothesis of $(2)$. \QEDA

\smallskip

With $X \subset \mathbb{P}^N$ and $Z \subset \mathbb{P}^{N-1}$ a hyperplane section of $X$, the following lemma gives a sufficient condition for the component of the Hilbert scheme containing $X$ to have a dominant map to the component of the Hilbert scheme containing $Z$.

\begin{lemma}\label{hyperplane sections fill out entire moduli}
Let $X \subset \mathbb{P}^{N}$ be a smooth projective variety and $Z \subset \mathbb{P}^{N-1}$ denote a smooth hyperplane section. If $X$ is unobstructed inside $\mathbb{P}^{N}$ and $H^1(N_{X/\mathbb{P}^{N}}(-1)) = 0$, then $Z$ is unobstructed inside $\mathbb{P}^{N-1}$ and there is a dominant map from the unique component of the Hilbert scheme in $\mathbb{P}^{N}$ that contains $X$ to the unique component of the Hilbert scheme of $\mathbb{P}^{N-1}$ that contains $Z$ which sends $X$ to $Z$.
\end{lemma}

\noindent\textit{Proof.} Let the unique component of the Hilbert scheme containing $X$ be denoted by $\mathcal{H}_X$. Let $H$ be the hyperplane such that $X \cap H = Z$. Now consider the universal family $\mathcal{F}_X \subset \mathbb{P}^N \times \mathcal{H}_X$ over $\mathcal{H}_X$ and look at $\mathcal{F}_X \cap (H \times \mathcal{H}_X)$, the intersection being the scheme theoretic intersection inside $\mathbb{P}^N \times \mathcal{H}_X$. This intersection gives a flat family of subschemes with the same Hilbert polynomial as $Z \subset \mathbb{P}^{N-1}$ parameterized by an open set $\mathcal{U}_X$ of $\mathcal{H}_X$. Hence we get a map from $\mathcal{U}_X$ to the unique component of the Hilbert scheme $\mathcal{H}_Z$ containing $Z$. We show that this map is dominant. For that it is enough to show that the map at the level of tangent spaces is surjective and the map at the level of obstruction spaces is injective (see \cite[Proposition $2.2.5$ (iii), Proposition $2.3.6$]{Ser06}). \textcolor{black}{This simultaneously shows that $Z$ is unobstructed in $\mathbb{P}^{N-1}$ using the fact that $X$ is unobstructed in $\mathbb{P}^N$ (see \cite[Proposition $2.2.5$ (iii), Proposition $2.3.6$]{Ser06}) } We have the following exact sequence:

$$0 \to N_{X/\mathbb{P}^{N}}(-1) \to N_{X/\mathbb{P}^{N}} \to N_{Z/\mathbb{P}^{N-1}} \to 0$$
Taking cohomology we have that 
$$H^0(N_{X/\mathbb{P}^{N}}) \to H^0(N_{Z/\mathbb{P}^{N-1}}) \to H^1(N_{X/\mathbb{P}^{N}}(-1)) \to H^1(N_{X/\mathbb{P}^{N}}) \to H^1(N_{Z/\mathbb{P}^{N-1}})$$
Now the vanishing of $H^1(N_{X/\mathbb{P}^{N}}(-1))$ proves our result since the leftmost map is the map of tangent spaces and the rightmost map is the map at the level of obstruction spaces. \QEDA

\smallskip

With $X \subset \mathbb{P}^N$ and $Z \subset \mathbb{P}^{N-1}$ a hyperplane section of $X$, the following lemma gives a sufficient condition for smoothness of the Hilbert point of the projective cone $C(Z) \subset \mathbb{P}^N$. 

\begin{lemma}\label{smoothness of Hilbert point of the cone}
   Suppose that $X \subset \mathbb{P}^{N}$ be a smooth projective variety and $Z \subset \mathbb{P}^{N-1}$ be a general hyperplane section which is assumed to be projectively normal. Assume
   \begin{itemize}
       \item[(1)] $X$ is unobstructed inside $\mathbb{P}^{N}$
       \item[(2)] $\alpha(X) = 0$ and $\alpha(Z) = 1$
       \item[(3)] $H^1(N_{X/\mathbb{P}^{N}}(-k)) = 0$ for $k = 1,2,3$.
   \end{itemize}
   Then the projective cone $C(Z) \subset \mathbb{P}^{N}$ is a smooth point of both the Hilbert scheme containing $C(Z)$ and the fiber over $Z$ of the the flag Hilbert scheme containing the pair $(X \subset \mathbb{P}^{N}, Z \subset \mathbb{P}^{N-1})$.   
\end{lemma}

\noindent\textit{Proof.} The dimension of the tangent space to the Hilbert scheme at the point $C(Z)$ is given by $\Sigma_{k \geq 0}  h^0(N_{Z/\mathbb{P}^{N-1}})(-k)$. We have the exact sequence
$$0 \to N_{X/\mathbb{P}^{N}}(-1) \to N_{X/\mathbb{P}^{N}} \to N_{Z/\mathbb{P}^{N-1}} \to 0$$
Since $\alpha(X) = 0$, we have that $h^0(N_{X/\mathbb{P}^{N}}(-2)) = 0$ (see \cite[Proof of Theorem $1.3$, Claim $2.3$]{Lop23}). Combined with the fact that $h^1(N_{X/\mathbb{P}^{N}}(-3)) = 0$, we have that $h^0(N_{Z/\mathbb{P}^{N-1}})(-2) = 0$. Now $h^0(N_{Z/\mathbb{P}^{N-1}})(-1)-N = \alpha(Z)$. Hence the dimension of the tangent space to the Hilbert scheme at the point $C(Z)$ is given by $h^0(N_{Z/\mathbb{P}^{N-1}})+\alpha(Z)+N = h^0(N_{Z/\mathbb{P}^{N-1}})+\alpha(X)+N+1 = h^0(N_{Z/\mathbb{P}^{N-1}})+h^0(N_{X/\mathbb{P}^{N}}(-1))$.

Now once again from the previous exact sequence, using $h^1(N_{X/\mathbb{P}^{N}}(-1)) = 0$,  we have that $h^0(N_{X/\mathbb{P}^{N}}(-1)) = h^0(N_{X/\mathbb{P}^{N}})-h^0(N_{Z/\mathbb{P}^{N-1}})$. Then the dimension of the tangent space to the Hilbert scheme at the point $C(Z)$ is $h^0(N_{X/\mathbb{P}^{N}})$ which is the dimension of the Hilbert scheme since $X$ is unobstructed inside $\mathbb{P}^{N}$. \par

The dimension of the tangent space to the fiber of the flag Hilbert scheme at the point $(C(Z), Z)$ is given by $\Sigma_{k \geq 1}  h^0(N_{Z/\mathbb{P}^{N-1}})(-k)$. This dimension by the above computation is $h^0(N_{Z/\mathbb{P}^{N-1}})(-1) = h^0(N_{X/\mathbb{P}^{N}})(-1)$ using the fact that $h^0(N_{X/\mathbb{P}^{N}}(-2)) = h^1(N_{X/\mathbb{P}^{N}}(-2)) = 0$. But now $H^1(N_{X/\mathbb{P}^{N}}(-1))$ is the obstruction space to deformations of $X$ inside $\mathbb{P}^{N}$ keeping $Z$ fixed, or in other words the vanishing of $H^1(N_{X/\mathbb{P}^{N}}(-1))$ implies that $(X,Z)$ is a smooth point of the fiber over $Z$ of the flag Hilbert scheme with dimension of tangent space at the point $(X,Z)$ being $h^0(N_{X/\mathbb{P}^{N}}(-1))$. This shows that the point $(C(Z), Z)$ is also a smooth point of the fiber. \QEDA

For lemmas \ref{ZL for surface sections}, \ref{hyperplane sections fill out entire moduli}, \ref{smoothness of Hilbert point of the cone}, to hold true one requires the condition $H^1(N_{X/\mathbb{P}^N}(-k)) = 0$ for $k = 1,2,3$. To show this vanishing for a general $X$ in an irreducible component of a Hilbert scheme, we will degenerate $X$ to a suitable ribbon $\widetilde{Y}$ inside $\mathbb{P}^N$ and prove the vanishing of $H^1(N_{\widetilde{Y}/\mathbb{P}^N}(-k))$. The following theorem shows how to reduce the vanishing of $H^1(N_{\widetilde{Y}/\mathbb{P}^N}(-k))$ to vanishings of the cohomology of certain twists of the tangent bundle on the reduced part $Y$ of the ribbon $\widetilde{Y}$.

\begin{proposition}\label{cohomology of twists of normal bundles}

Let $\widetilde{Y}$ be a ribbon over a smooth projective variety $Y$ of dimension $d$ with conormal bundle $L$. Let $\widetilde{Y} \hookrightarrow \mathbb{P}^N$ be an embedding of $\widetilde{Y}$ induced by the complete linear series of a very ample line bundle $\widetilde{H}$. Let $H = \widetilde{H}|_Y$. If the following conditions are satisfied 
\begin{itemize}
    \item[(1)] $H^1(\mathcal{O}_Y(-L-kH)) = H^1(\mathcal{O}_Y(-2L-kH)) = H^2(\mathcal{O}_Y(-kH)) = H^2(\mathcal{O}_Y(-L-kH)) = \\ H^1(\mathcal{O}_Y(L-(k-1)H) = H^1(\mathcal{O}_Y(-(k-1)H) = H^2(\mathcal{O}_Y(L-kH) = H^2(\mathcal{O}_Y(-kH) = 0$
    \item[(2)] $H^2(T_Y(L-kH)) = H^2(T_Y(-kH)) = 0$ 
\end{itemize}
    
then, $H^1(N_{\widetilde{Y}/\mathbb{P}^N}(-k\widetilde{H})) = 0$ for $1 \leq k \leq 3$.

\end{proposition}

\noindent\textit{Proof.} \textcolor{black}{The idea of the proof is to reduce the needed vanishing on $\widetilde Y$ to relevant cohomology  vanishings on $Y$. To carry this out, we use the following exact sequences}:
\begin{equation}\label{a}
    0 \to N_{\widetilde{Y}/\mathbb{P}^N} (L) \to N_{\widetilde{Y}/\mathbb{P}^N} \to N_{\widetilde{Y}/\mathbb{P}^N} \otimes \mathcal{O}_Y \to 0
\end{equation}
\begin{equation}\label{b}
    0 \to \textrm{Hom}(I_{\widetilde{Y}}/I_{Y}^2, \mathcal{O}_Y) \to N_{\widetilde{Y}/\mathbb{P}^N} \otimes \mathcal{O}_Y \to \mathcal{O}_Y(-2L) \to 0 
\end{equation}
\begin{equation}\label{c}
    0 \to \mathcal{O}_Y(-L) \to N_{Y/\mathbb{P}^N} \to \textrm{Hom}(I_{\widetilde{Y}}/I_{Y}^2, \mathcal{O}_Y) \to 0
\end{equation}
It is enough to show the vanishings of $H^1(N_{\widetilde{Y}/\mathbb{P}^N}(L-kH)) = H^1(N_{\widetilde{Y}/\mathbb{P}^N}(-kH)) = 0$ for $1 \leq k \leq 3$.  We show how to derive the conditions for the first vanishing. From the exact sequence \ref{b}, it is enough to show  $H^1(\textrm{Hom}(I_{\widetilde{Y}}/I_{Y}^2, \mathcal{O}_Y)(L-kH)) = H^1(\mathcal{O}_Y(-L-kH)) = 0$. From the exact sequence \ref{c}, to show the first vanishing, it is enough to show $H^1(N_{Y/\mathbb{P}^N}(L-kH)) = H^2(\mathcal{O}_Y(-kH)) = 0$. To show the first of the vanishings, using the exact sequence

$$0 \to T_Y \to T_{\mathbb{P}^N|_Y} \to N_{Y/\mathbb{P}^N} \to 0$$
it is enough to show $H^1(T_{\mathbb{P}^N|_Y}(L-kH)) = H^2(T_Y(L-kH)) = 0$.
To show the first of the vanishings, using the Euler sequence it is enough to show $H^1(\mathcal{O}_Y(L-(k-1)H) = H^2(\mathcal{O}_Y(L-kH) = 0$. \textcolor{black}{The second vanishing, namely $H^1(N_{\widetilde{Y}/\mathbb{P}^N}(-kH)) = 0$ follows from a similar procedure as above, and we leave it to the reader.} This completes the proof. \QEDA


\section{Calabi-Yau ribbons on Fano threefolds and deformations of Calabi-Yau double covers of Fano threefolds}

\textcolor{black}{A ribbon} $\widetilde{Y}$ is a Calabi-Yau ribbon on $Y$ if and only if $L = \mathcal{O}_Y(K_Y)$. Let $Y$ be a smooth, projective Fano threefold of Picard rank $1$. Let $H$ denote the generator of Pic $(Y)$. Let $j$ be the least positive integer such that $j$ is greater than or equal to the index of $Y$ and $jH$ is very ample. \textcolor{black}{Let $i(Y)$ denote the index of the fano threefold $Y$.} In the following proposition we show that for each deformation type of a Fano threefold $Y$ and a positive integer $l \geq j$, there exist a unique component $\mathcal{H}^Y_l$ of the Hilbert scheme parameterizing smooth Calabi-Yau threefolds in general and contains Calabi-Yau ribbons as a special locus.


    

\begin{proposition}\label{existence and smoothing of ribbons}
 Let $Y$ be a smooth, projective Fano threefold of Picard rank $1$. Let $\pi: X \to Y$ denote a Calabi-Yau double cover branched along a smooth divisor $B \in |-2K_Y|$. Let $H$ denote the generator of Pic $(Y)$. \textcolor{black}{Let $j \geq i(Y)$ be the least positive integer such that} $jH$ is very ample. Then
\begin{enumerate}
     \item Pic $(X) = \mathbb{Z}$. If $A$ denotes the generator of Pic $(X)$, then for any polarized deformation $(X_t, A_t)$ of $(X,A)$ along a one parameter family $T$ the line bundle $A_t$ generates Pic $(X_t)$. 
     \item For $l \geq j$, $h^0(lA_t) = h^0(lH)+h^0((l-i(Y))H)$. We set $N_l+1 = h^0(lH)+h^0((l-i(Y)H)$.
     \item \textcolor{black}{Let $l \geq j$. If $Y$ belongs to family $1.11$, further assume $l \geq 4$. Then for a general polarized deformation $(X_t, A_t)$ of $(X,A)$ along a one parameter family $T$, $lA_t$ is very ample. Further} there exists a unique irreducible component $\mathcal{H}^Y_l$ of the Hilbert scheme of $\mathbb{P}^{N_l}$ with Hilbert polynomial $p(z) = h^0(lzH)+h^0((lz-i(Y))H)$ that parameterizes linearly normal Calabi-Yau varieties $X_t \subset \mathbb{P}^{N_l}$ such that $(X_t, \mathcal{O}_{X_t}(1) = lA_t)$ can be deformed to $(X, lA)$ along an irreducible one parameter family.
     \item The component $\mathcal{H}^Y_l$ contains linearly normal embedded Calabi-Yau ribbons $\widetilde{Y} \subset \mathbb{P}^{N_l}$ on $Y$ such that $\mathcal{O}_{\widetilde{Y}}(1)|_Y = lH$. \textcolor{black}{The embedded ribbons form an irreducible family, and they represent smooth points} of $\mathcal{H}^Y_l$. \textcolor{black}{Consequently, any linearly normal embedded Calabi-Yau ribbon $\widetilde{Y} \subset \mathbb{P}^{N_l}$ on $Y$ such that $\mathcal{O}_{\widetilde{Y}}(1)|_Y = lH$ is contained in $\mathcal{H}^Y_l$}

\end{enumerate}

\end{proposition} 


\textcolor{black}{The following lemma shows optimality of the above proposition, and gives precise numbers on the very ampleness of $lA_t$ unless $Y$ is of type $1.1$ and $l = 2$. When $Y$ is of type $1.1$ and $l = 2$, then $2A_t$ is very ample on $X_t$ but there does not exist a Calabi-Yau ribbon supported on $Y$ in $\mathcal{H}^Y_2$.}


\begin{lemma}\label{values of j}

\textcolor{black}{Let the notation be as in Proposition 3.1. Assume if $Y$ is of type $1.1$, then $l \neq 2$. Then $lA_t$ is very ample on $X_t$ if and only if $l \geq j$. The $j$ are as follows for various families.}
\begin{enumerate}
    \item For family $1.1$, $j = 3$. 
    \item For families $1.2-1.10$, $j = 1$. For families $1.2-1.4$, the general elements of $\mathcal{H}^Y_1$ are complete intersections.
    \item For family $1.11$, $j = 3$. 
    \item For family $1.12-1.15$, $j = 2$. 
    \item For family $1.16$, $j = 3$. 
    \item For family $1.17$, $j = 4$. 
\end{enumerate}  
\textcolor{black}{Moreover, whenever $lA_t$ is very ample the embedding of $X_t$ by $|lA_t|$, degenerates to a Calabi-Yau ribbon on $Y$.}

\end{lemma}

\noindent\textit{Proof.} Let us first show the values of $j$. For families $1.2-1.10$, $1.12-1.17$, this is clear by \cite{Bel24} since $-K_Y$ is very ample. For $1.1$, $H = -K_Y$. $Y$ is a double cover of $\mathbb{P}^3$ branched along a sextic. We have $h^0(-2K_Y) = h^0(\mathcal{O}_{\mathbb{P}^3}(2))$ and hence the morphism induced by $|-2K_Y|$ is $2:1$ onto the second Veronese embedding of $\mathbb{P}^3$. Further one can check that $-3K_Y$ is in fact projectively normal and hence $j = 3$. For $1.11$, by \cite{Bel24}, $j \geq 3$. Now by \cite[Theorem $0.2$]{NO02}, $lH$ is very ample on $Y$ for $l \geq 3$, which implies $j = 3$. 

\par

Now we show that if $l < j$, then $lA_t$ is not very ample. \textcolor{black}{For family $1.1$, $l = 1$, if $A_t$ is very ample, then it embeds $X_t$ as a degree $4$ Calabi-Yau threefold inside $\mathbb{P}^4$. This gives a contradiction.} \textcolor{black}{For family $1.11$, $l = 1$, $A_t$ is not very ample since $h^0(A_t) = 3$. For $l = 2$, if $|2A_t|$ is very ample, then $X_t$ is an aCM Calabi-Yau threefold of degree $16$ inside $\mathbb{P}^6$. But such a Calabi-Yau does not exist due to
\cite[Theorem $3.2$]{KK16}.} For family $1.12-1.15$, $j = 2$. For $l = 1$, $A_t$ is not very ample. \textcolor{black}{To see this first note that $h^0(A_t) = h^0(H)$. Now for $1.12$, $h^0(A_t) = 4$ and the morphism given by $|A_t|$ is $4:1$ onto $\mathbb{P}^3$. For $1.13-1.14$, if $A_t$ is very ample, then $|A_t|$ embeds $X_t$ inside $\mathbb{P}^5$. By \cite[Remark $11$]{BC83}, all Calabi-Yau threefolds inside $\mathbb{P}^5$ are complete intersections of either two cubics or a quadric and a quartic or a quintic and a hyperplane. For $1.13$, $A_t^3 = 6$, which cannot be the degree of a Calabi-Yau inside $\mathbb{P}^5$. For $1.14$, $A_t^3 = 8$, so $X_t$ has to be the complete intersection of a quadric and a quartic. But then this is not possible by comparing $h^{1,2}(X_t)$ with $h^{1,2}$ of the complete intersection. For $1.15$, if $A_t$ is very ample then $|A_t|$ embeds $X_t$ as a degree $10$, aCM Calabi-Yau threefold inside $\mathbb{P}^6$. By \cite[Theorem $3.2$, Section $6$]{KK16}, such a Calabi-Yau threefold cannot exist.} For family $1.16$, $j = 3$. For $l = 1,2$, $lA_t$ is not very ample. \textcolor{black}{Note first that $H^1(T_Y \otimes K_Y) = H^{1,2}(Y) = 0$. Now since $h^0(lA) = h^0(lH)$, the morphism given by $|A|$ on the double cover $X$ factors through a $2:1$ map onto $Y$ followed by $Y$ embedded by $|lH|$. This together with the fact that $H^1(\mathcal{O}_{X}) = 0$, implies by \cite[Corollary $1.1$]{Weh86} and arguments similar to \cite[Theorem $1.7$]{GGP16} the morphism by $|lA_t|$ factors as a double cover of $Y$ as well.} For family $1.17$, $j = 4$. \textcolor{black}{For $l = 1,2,3$, the map given by $|A_t|$ is not very ample by exactly same arguments as in part $(5)$.} The last statement follows by combining Lemma \ref{values of j} and Proposition \ref{existence and smoothing of ribbons}.  \QEDA

\vspace{0.3cm}

\noindent\textit{Proof of Proposition \ref{existence and smoothing of ribbons}.} 
We first prove part $(1)$. Consider the Calabi-Yau double cover $\pi: X \to Y$ branched along $B \in |-2K_Y|$. By \cite{Par91}, Lemma $4.2$, we have $H^1(\Omega_X^1) = H^1(\pi_*\Omega_X^1) = H^1(\Omega_Y) \bigoplus H^1(\Omega_Y(\log (B)) \otimes K_Y)$. Now using the exact sequence 
$$ 0 \to \Omega_Y \otimes K_Y \to \Omega_Y(\log (B)) \otimes K_Y \to K_Y|_B \to 0 $$ we have 
$$H^1(\Omega_Y \otimes K_Y) \to H^1(\Omega_Y(\log (B)) \otimes K_Y) \to H^1(K_Y|_B)$$
The leftmost term is zero due to Nakano vanishing theorem. We have 
$$0 \to 3K_Y \to K_Y \to K_Y|_B \to 0$$ 
and hence we see that $H^1(K_Y|_B) = 0$. This implies $H^1(\Omega_Y(\log (B)) \otimes K_Y) = 0$. Now since $H^1(\mathcal{O}_Y) = H^1(\mathcal{O}_X) = 0$, the pullback map between $\pi^*: \textrm{Pic}(Y)  \to \textrm{Pic}(X) $ is induced by the map 
$$H^1(\Omega_Y) \to H^1(\Omega_Y) \bigoplus H^1(\Omega_Y(\log (B)) \otimes K_Y)$$
Therefore $H^1(\Omega_Y(\log (B)) \otimes K_Y) = 0$ implies $\pi^*$ is an isomorphism, $h^{1,1}(X) = 1$ and $\textrm{Pic}(X) = \mathbb{Z}$. Let $A$ denote the generator. We have that for $X_t$, a general deformation of $X$, $h^{1,1}(X_t) = 1$. On the other hand, by \cite[Corollary 4.18]{Debarre} and \cite[Theorem 0.1]{KMM},  $Y$ is simply connected, so, by the Lefschetz theorem  for homotopy groups (see \cite[Theorem 3.1.21]{positivity}), so is $B$, where $B$ is the ample branch divisor of the double cover $X$. Now if $R \subset X$ denotes the ramification divisor, $R \cong B$ and hence $R$ is simply connected. Since $R$ is ample, once again by the Lefschetz theorem  for homotopy groups, so is $X$. Therefore $X$ is simply connected and hence so is $X_t$. Therefore, $\textrm{Pic}(X_t) = \mathbb{Z}$. Again since divisibility of a line bundle remains constant along a smooth family, we have that for any polarized deformation $(X_t, A_t)$ of $(X,A)$, $A_t$ generates the Picard group of $X_t$. \par
\textcolor{black}{We now prove part $(2)$. We have by part $(1)$, $h^0(lA) = h^0(\pi^*(lH)) = h^0(lH) + h^0((lH) \otimes K_Y) = h^0(lH)+h^0((l-i(Y))H)$. Now since by Kodaira vanishing Theorem, $H^1(lA) = 0$, we have that for a general polarized deformation $(X_t,A_t)$ of $(X,A)$, along a one -parameter family $T$, $h^0(lA_t) = h^0(lH)+h^0((l-i(Y))H)$.} \par
\textcolor{black}{We now prove part $(3)$ and $(4)$. To do that,}
we first show that there is a unique irreducible component $\mathcal{H}^Y_l$ of the Hilbert scheme inside $\mathbb{P}^{N_l}$ which contains all linearly normal Calabi-Yau ribbons $\widetilde{Y} \subset \mathbb{P}^{N_l}$ such that $\mathcal{O}_{\widetilde{Y}}(1)|_Y = lH$. Let 
$Y \subset \mathbb{P}^{N}$ be the embedding given by the complete linear series of $lH$. \textcolor{black}{By \cite[Proposition $2.1$]{Gon06}}, Calabi-Yau ribbons $\widetilde{Y} \subset \mathbb{P}^N$ extending the embedding $Y \subset \mathbb{P}^N$ are parameterized by nowhere vanishing sections of $H^0(N_{Y/\mathbb{P}^N} \otimes K_Y)$. We have the exact sequence 
$$0 \to K_Y \to (l-i(Y))H^{N+1} \to T_{\mathbb{P}^N}|_Y \otimes K_Y \to 0$$
\textcolor{black}{Since we have $l \geq j$ and if $Y$ belongs to family $1.11$, $l \geq 4$, we have using the values of $j$ from \Cref{values of j} and \cite{Bel24}, that $(l-i(Y))H$ is base point free. } \textcolor{black}{Hence} $T_{\mathbb{P}^N}|_Y \otimes K_Y$ is base point free. Now from the exact sequence 
$$0 \to T_Y \otimes K_Y \to T_{\mathbb{P}^N}|_Y \otimes K_Y \to N_{Y/\mathbb{P}^N} \otimes K_Y$$
we conclude that $N_{Y/\mathbb{P}^N} \otimes K_Y$ is base point free. Now the rank of $N_{Y/\mathbb{P}^N} \otimes K_Y$ which is the codimension of $Y$ inside $\mathbb{P}^N$ is greater than the dimension of $Y$, which is $3$. Therefore the vector bundle must have a no-where vanishing section. This corresponds to a ribbon $\widetilde{Y} \subset \mathbb{P}^N$ \textcolor{black}{with conormal bundle $K_Y$,} extending the embedding $Y \subset \mathbb{P}^N$. 
The embedding $\widetilde{Y} \subset \mathbb{P}^N$ is induced by some sublinear series of a very ample line bundle $\widetilde{H}$ on $\widetilde{Y}$ such that $\widetilde{H}|_Y = lH$. Re-embedding $\widetilde{Y}$ be the complete linear series of $\widetilde{H}$ gives an embedding $\widetilde{Y} \subset \mathbb{P}^{N'}$ (such that the reduced part $Y$ is now embedded degenerately inside a $\mathbb{P}^N \subset \mathbb{P}^{N'})$. \textcolor{black}{All ribbons $\widetilde{Y} \subset \mathbb{P}^{N'}$ on $Y$, with conormal bundle $K_Y$, extending the degenerate embedding $Y \subset \mathbb{P}^{N'}$ forms an irreducible family in $\mathbb{P}^{N'}$. This is due to the fact that all such ribbons $\widetilde{Y} \subset \mathbb{P}^{N'}$ are parameterized by an open subset of $\mathbb{P}(H^0(N_{Y/\mathbb{P}^{N'}} \otimes K_Y)$, given by no-where vanishing sections.} Let us show $N' = N_l = h^0(lH)+h^0((l-i(Y))H)$. We have an exact sequence  $$0 \to K_Y \to \mathcal{O}_{\widetilde{Y}} \to \mathcal{O}_Y \to 0$$ 
which after tensoring with $\widetilde{H}$ and taking cohomology yields
$$0 \to H^0(lH \otimes K_Y) \to H^0(\widetilde{H}) \to H^0(lH) \to 0$$ 
Noting that $H^0(lH \otimes K_Y) = H^0((l-i(Y))H)$ and $H^1((l-i(Y))H) = 0$, we have $N'+1 = h^0(\widetilde{H}) = h^0(lH)+h^0((l-i(Y))H) = N_l+1$.
Now we show that all such ribbons $\widetilde{Y} \subset \mathbb{P}^{N_l}$ represent smooth points of the Hilbert scheme. This can be checked using \cite[Lemma $2.14$]{BGM24}. \color{black} The only non-trivial vanishings are those of $H^1(N_{Y/\mathbb{P}^{N_l}})$ and $H^1(N_{Y/\mathbb{P}^{N_l}} \otimes K_Y$). We use the exact sequence 

$$0 \to T_Y \to T_{\mathbb{P}^{N_l}|_Y} \to N_{Y/\mathbb{P}^{N_l}} \to 0$$
Using Kodaira vanishing theorem, and the Euler exact sequence of $\mathbb{P}^{N_l}$ restricted to $Y$, we have that $H^1(T_{\mathbb{P}^{N_l}|_Y}) = H^1(T_{\mathbb{P}^{N_l}|_Y} \otimes L) = 0$. Now note that by Nakano-vanishing theorem, $H^2(T_Y) = H^1(\Omega_Y \otimes K_Y)^* = 0$. This implies that $H^1(N_{Y/\mathbb{P}^{N_l}}) = 0$. To show the vanishing of $H^1(N_{Y/\mathbb{P}^{N_l}} \otimes K_Y)$, first note that the vanishings so far yield an exact sequence 
$$0 \to H^1(N_{Y/\mathbb{P}^{N_l}} \otimes K_Y) \to H^2(T_Y \otimes K_Y) \xrightarrow{f} H^2(T_{\mathbb{P}^{N_l}|_Y} \otimes K_Y)$$
We show that $f$ is injective. 
The homomorphism $f$ can be identified, via Serre duality, with the dual of the homomorphism
\begin{equation*}\label{$(1)$}
    H^1(\Omega_{\mathbb{P}^{N_l}}|_Y) \overset{g} {\longrightarrow} H^1(\Omega_Y).
\end{equation*}
Composing $g$ with the natural map 
$H^1(\Omega_{\mathbb{P}^{N_l}}) \to H^1(\Omega_{\mathbb{P}^{N_l}}|_Y)$, we get 
a homomorphism 
\begin{equation*}
    H^1(\Omega_{\mathbb{P}^{N_l}}) \overset{\widehat g} {\longrightarrow} H^1(\Omega_Y).
\end{equation*}
which induces, by restriction,  a homomorphism between the Neron-Severi groups of $\mathbb P^{N_l}$ and $Y$, which, in our situation, is a homomorphism $h$ between the Picard groups of $\mathbb P^{N_l}$ and $Y$, because $Y$, being Fano, is a regular variety. \textcolor{black}{The homomorphism $h$ sends a line bundle on $\mathbb P^{N_l}$ to its restriction on $Y$}, so it is nonzero. Then, so is $\widehat g$ and, consequently, so is $g$. Therefore, $f$ is also nonzero. Since, by Serre duality, $h^2(T_Y \otimes \omega_Y) = h^1(\Omega_Y)$ and, by assumption, $h^{1,1}(Y) = 1$, the homomorphism $f$ is injective as desired (in fact, it is an isomorphism). 
This finishes the proof of the smoothness of the Hilbert point of $\widetilde{Y}$. \par
\color{black}
\textcolor{black}{Thus so far we have shown the existence of a unique component $\mathcal{H}^Y_l$ of the Hilbert scheme of $\mathbb{P}^{N_l}$, with Hilbert polynomial $p(z) = h^0(lzH)+h^0((lz-i(Y))H)$, containing any linearly normal embedded Calabi-Yau ribbon $\widetilde{Y} \subset \mathbb{P}^{N_l}$ on $Y$ such that $\mathcal{O}_{\widetilde{Y}}(1)|_Y = lH$.} \par
Now we show that for $t \in \mathcal{H}^Y_l$ general, the varieties $X_t \subset \mathbb{P}^{N_l}$ are smooth linearly normal Calabi-Yau varieties such that $(X_t, \mathcal{O}_{X_t}(1))$ is a polarized deformation $(X,lA)$. As we have shown before, the embedding $\widetilde{Y} \subset \mathbb{P}^{N_l}$ induces a degenerate embedding $Y \subset \mathbb{P}^{N_l}$. Composing with the map $\pi: X \to Y$ we have a diagram

\[
\begin{tikzcd}
X \arrow[dr, "\varphi"] \arrow[d, "\pi"] & \\
Y \arrow[r, hook] & \mathbb{P}^{N_l}
\end{tikzcd}
\]

We show that the morphism $\varphi$ can be deformed to an embedding $\varphi_t$ along a one-parameter family $t \in T$. To see this note that this diagram has an associated exact sequence 
$$0 \to N_{\pi} \to N_{\varphi} \to \pi^*N_{Y/\mathbb{P}^{N_l}} \to 0$$
Taking cohomology and noting that $H^1(N_{\pi}) = H^1(B|_B) = 0$ \textcolor{black}{(see for example \cite[Lemma $2.5$]{GGP10})}, we have a surjection
$$H^0(N_{\varphi}) \to H^0(N_{Y/\mathbb{P}^{N_l}}) \oplus H^0(N_{Y/\mathbb{P}^{N_l}} \otimes K_Y) \to 0$$

So one can pick a first order deformation of $\varphi$ that maps to a surjective homomorphism in $H^0(N_{Y/\mathbb{P}^{N_l}} \otimes K_Y)$ which corresponds to a class of ribbon $\widetilde{Y} \subset \mathbb{P}^{N_l}$. Further we have

$$0 \to H^1(N_{\pi}) \to H^1(N_{\varphi}) \to H^1(N_{Y/\mathbb{P}^{N_l}}) \oplus H^1(N_{Y/\mathbb{P}^{N_l}} \otimes K_Y)$$

As before since $Y$ is Fano and using Nakano-Vanishing theorem, one can check that both flanking terms of the previous exact sequence vanishes. This implies that $H^1(N_{\varphi}) = 0$ and deformations of $\varphi$ are unobstructed.
This implies by \cite[Proposition $1.4$]{GGP13}, that there is a one parameter family $\Phi_T: \mathcal{X} \to \mathbb{P}^{N_l}_T$ such that $\Phi_0 = \varphi$ while $\Phi_t = \varphi_t: \mathcal{X}_t = X_t \to \mathbb{P}^{N_l}$ is an embedding. Further (Im $\Phi)_0 = \widetilde{Y} \subset \mathbb{P}^{N_l}$ while for $t \neq 0$, (Im $\Phi)_t = X_t \subset \mathbb{P}^{N_l}$. Note that since $\widetilde{Y} \subset \mathbb{P}^{N_l}$ is linearly normal, we have that
$X_t \subset \mathbb{P}^{N_l}$ is also linearly normal. Further we have $(\Phi_T^*\mathcal{O}_{\mathbb{P}^{N_l}_T}(1))_0 = \varphi^*\mathcal{O}_{\mathbb{P}^{N_l}}(1) = \pi^*(lH) = lA$, and hence $(\Phi_T^*\mathcal{O}_{\mathbb{P}^{N_l}_T}(1))_t = \mathcal{O}_{X_t}(1) = lA_t$ is deformation of $A$. \par
\textcolor{black}{We end the proof by showing that if there is a Calabi-Yau threefold $X_t \subset \mathbb{P}^{N_l}$ such that the pair $(X_t, \mathcal{O}_{X_t}(1))$ can be deformed to $(X,lA)$ along an irredicible curve $T$, then $X_t \subset \mathbb{P}^{N_l}$ belongs to the irreducible component $\mathcal{H}_l^Y$. } Suppose that there are two such irreducible components. Pick $X_t \subset \mathbb{P}^{N_l}$ from one component and $X'_s \subset \mathbb{P}^{N_l}$ from the other. Now consider the deformations of the polarized pair $(X,lA)$. We have the Atiyah extension of $lA$ given by

$$0 \to \mathcal{O}_X \to \mathcal{E}_{lA} \to T_X \to 0$$

Taking cohomology of the sequence we have that 

$$H^1(\mathcal{E}_{lA}) \to H^1(T_X) \to H^2(\mathcal{O}_X) \to H^2(\mathcal{E}_{lA}) \to H^2(T_X)$$

Since $H^2(\mathcal{O}_X) = 0$, we have that the map $H^1(\mathcal{E}_L) \to H^1(T_X)$ between the first order deformations of the polarized pair $(X,lA)$
and the first order deformations of $X$ is surjective while the respective maps of obstruction spaces is injective. So we have that the forgetful map between the functor of deformations of the polarized pair to the functor of deformations of $X$ is smooth.
By the $T^1-$ lifting criterion by \cite{Ran92}, $X$ is unobstructed and hence $(X,lA)$ is unobstructed as well. So the pair $(X, lA)$ is contained in a unique irreducible component of the moduli space of smooth polarized Calabi-Yau threefolds. Since both of the pairs $(X_t, \mathcal{O}_{X_t}(1))$ and $(X'_s, \mathcal{O}_{X'_s}(1))$ can be deformed to $(X,lA)$ along an irreducble curve, both of them are  contained the same component of smooth polarized Calabi-Yau threefolds. So there is a one parameter family $(\mathcal{X}, \mathcal{O}_{\mathcal{X}}(1))$ where $f: \mathcal{X} \to T$ is a scheme over $T$, such that $(\mathcal{X}_t, \mathcal{O}_{\mathcal{X}_t}(1)) = (X_t, \mathcal{O}_{X_t}(1))$ while $(\mathcal{X}_s, \mathcal{O}_{\mathcal{X}_s}(1)) = (X'_s, \mathcal{O}_{X'_s}(1))$. But now since $H^i(\mathcal{O}_{\mathcal{X}_p}(1))$ for $i > 0$ for any $p \in T$, we have that $h^0(\mathcal{O}_{\mathcal{X}_p}(1))$ is constant for $p \in T$. Therefore $f_*(\mathcal{O}_{\mathcal{X}}(1))$ is locally free of rank $N_l+1$ and $\mathcal{X} \to \mathbb{P}(f_*(\mathcal{O}_{\mathcal{X}}(1)))$ is a $T-$ morphism for $T$ irreducible such that at $t \in T$ we get the embedding $X_t \subset \mathbb{P}^{N_l}$ while at $s \in T$, we get the embedding $X'_s \subset \mathbb{P}^{N_l}$. Therefore the $X_t \subset \mathbb{P}^{N_l}$ and $X'_s \subset \mathbb{P}^{N_l}$ belong to the same irreducible component of the Hilbert scheme. \QEDA

\section{Applications to extendability of Calabi-Yau threefolds and canonical surfaces}

\smallskip

\subsection{Extendability of general elements of $\mathcal{H}^Y_l$ where $l = j$ and $j$ equals the index of $Y$}

When $j$ equals the index $i(Y)$ of $Y$, i.e, when $-K_Y$ is very ample, we use the classification of Mukai varieties, to show that the ribbons $\widetilde{Y}$ in $\mathcal{H}^Y_j$ in Proposition \ref{existence and smoothing of ribbons}, $(4)$, are at least as much extendable as the underlying Fano variety $Y$ and eventually conclude the same for the general smoothing of $\widetilde{Y}$. \textcolor{black}{The metaprinciple 2.1 will guide us in the proof of the following theorem.}





\begin{theorem}\label{l=j}
\textcolor{black}{When $j = i(Y)$, the general Calabi-Yau threefolds in $\mathcal{H}^Y_j$ are at least as many times smoothly extendable as the Fano threefold $Y$ in its anticanonical embedding.}
\textcolor{black}{Using} notation as in Proposition \ref{existence and smoothing of ribbons} \textcolor{black}{and the numbering of the families of Fano-threefolds as in \url{https://www.fanography.info/}}, the extendability of the general Calabi-Yau threefolds in $\mathcal{H}^Y_j$ is summarized in the following table   
\end{theorem}

\begin{center}\phantomsection\label{t01}
\begin{tabular}{c|c|c} 
\hline
\makecell{Deformation type \\ of Fano threefold \\ $Y$} & \makecell{Value of \\ $l = j$}  & \makecell{$k$-extendability of \\ a general Calabi-Yau \\ threefold $X_t$ \\ in $\mathcal{H}_j^Y$}   \\
\hline\hline

$1.2-1.4$ & $l = 1$  & \makecell{infinitely many times \\ smoothly extendable} \\

\hline

$1.5$ & $l = 1$  & \makecell{smoothly $3$- extendable} \\

\hline

$1.6$ & $l = 1$  & \makecell{smoothly $7$- extendable} \\

\hline 
$1.7$ & $l = 1$  & \makecell{smoothly $5$- extendable} \\

\hline

$1.8$ & $l = 1$  & \makecell{smoothly $3$- extendable} \\

\hline

$1.9$ & $l = 1$  & \makecell{smoothly $2$- extendable} \\

\hline



$1.14$ & $l = 2$ & \makecell{smoothly $2$- extendable} \\

\hline





\end{tabular}
\end{center}
\captionof{table}{}
\vspace{0.3cm}

\noindent\textit{Proof.} \textcolor{black}{First of all, note that by \cite{M89A}, \cite{M89B} or \cite{CLM98} or \cite{Bel24}, the only families of Fano varieties that are smoothly-extendable are the ones listed in the table.}  In the case, $1.2-1.4$, the general elements of $\mathcal{H}^Y_1$ are complete intersections and hence the result follows. We concentrate on $1.5-1.9$ and $1.14$.  We first show that in each case, a special locus in $\mathcal{H}^Y_j$ is $k-$ extendable for $k$ as mentioned in the table. In each of the cases the $k-$ th step of the respective embedding $Y \xhookrightarrow[]{|\mathcal{O}_{Y}(j)|} \mathbb{P}^M$ is a Fano variety $G \hookrightarrow \mathbb{P}^{M+k}$ as described in \cite{M89A}, \cite{M89B} or \cite{CLM98} or \cite{Bel24}.

\vspace{0.3cm}

If $\mathcal{O}_G(1)$ be the line bundle giving the embedding of $G$ and $\mathcal{O}_Y(1)$ denote its pullback to $Y$, then $K_Y = \mathcal{O}_Y(-1)$.
We first show that there exists polarized ribbons $(\widetilde{G}, \mathcal{O}_{\widetilde{G}}(1))$ on $G$ with conormal bundle $\mathcal{O}_G(-1)$ such that $\mathcal{O}_{\widetilde{G}}(1)$ is very ample and $\mathcal{O}_{\widetilde{G}}(1)|_G = \mathcal{O}_G(1)$. For this consider the embedding $G \subset \mathbb{P}^{M+k}$ induced by the complete linear series of $\mathcal{O}_G(1)$. We have the twisted exact sequence of tangent bundles

$$0 \to T_G(-1) \to T_{\mathbb{P}^{M+k}}|_G(-1) \to N_{G/\mathbb{P}^{M+k}}(-1) \to 0$$


Again by twisting the pullback of the Euler exact sequence to $G$ by $\mathcal{O}_G(-1)$ we have 


$$0 \to \mathcal{O}_G(-1) \to \mathcal{O}_G^{\oplus {M+k}} \to T_{\mathbb{P}^{M+k}|_G}(-1) \to 0 $$

This implies that $T_{\mathbb{P}^{M+k}|_G}(-1)$ is base point free and hence $N_{G/\mathbb{P}^{M+k}}(-1)$ is base point free. Further the rank of $N_{G/\mathbb{P}^{M+k}}(-1)$ which is $M-3$ is greater than $3$. Hence there exists a nowhere vanishing section of $H^0(N_{G/\mathbb{P}^{M+k}}(-1))$ which corresponds to an embedded ribbon $\widetilde{G} \subset \mathbb{P}^{M+k}$ on $G \subset \mathbb{P}^{M+k}$. Hence we get a very ample polarization $\mathcal{O}_{\widetilde{G}}(1)$ which restricts to $\mathcal{O}_G(1)$ on $G$. Note that $H^0(\mathcal{O}_{\widetilde{G}}(1)) = H^0(\mathcal{O}_G(1))+1$ and hence the embedding $\widetilde{G} \subset \mathbb{P}^{M+k}$ is actually induced by a sublinear series of $\mathcal{O}_{\widetilde{G}}(1)$ of codimension $1$. Now embed $\widetilde{G} \subset \mathbb{P}^{M+k+1}$ by the complete linear series of $\mathcal{O}_{\widetilde{G}}(1)$. To show that $\widetilde{G}$ has an embedded smoothing in $\mathbb{P}^{M+k+1}$, we work as in Proposition \ref{existence and smoothing of ribbons}. By \cite[Proposition $2.16$]{BGM24} (which uses \cite[Proposition $1.4$]{GGP13}) we need to show that $|\mathcal{O}_G(2)|$ has a smooth member and prove the vanishing of $H^1(\mathcal{O}_G(2)), H^1(N_{G/\mathbb{P}^M}(-1))$ and $H^1(N_{G/\mathbb{P}^M})$. \textcolor{black}{Since $\mathcal{O}_G(1)$ is very ample, the linear system of $|\mathcal{O}_G(2)|$ has a smooth member by Bertini's theorem.} The first one of the vanishings follows from Kodaira vanishing theorem. From the above two exact sequences, for the other two vanishings, one needs to show the vanishing of $H^1(T_{\mathbb{P}^M}|_G(-i))$ and $H^2(T_G(-i))$ for $i = 0,1$. The first one boils down to the vanishings of $H^1(\mathcal{O}_G(i))$ for $i = 0,1$ and $H^2(\mathcal{O}_G(-i))$ for $i = 0, 1$, both of which are true since $G$ is Fano. Finally $H^2(T_G(-i)) = H^{d-2}(\Omega_G \otimes K_G \otimes \mathcal{O}_G(-i)) = 0$ (where $d$ is the dimension of $G$) by Nakano vanishing theorem since $K_G \otimes \mathcal{O}_G(-i)$ is negative ample for $i = 0,1$. 

\par 

\vspace{0.3cm}

We now show that if some $X \in \mathcal{H}^Y_1$ is $k-$ extendable, then a general $X_t \in \mathcal{H}^Y_1$ element is \textcolor{black}{$k$-extendable}. \textcolor{black}{Let} $X_{k+1} \supset X_{k} \supset ... \supset X_1 = X$ be the chain of extensions. Since $X_1$ is a prime polarized Calabi-Yau threefold, $X_2$ is an anticanonically embedded Fano fourfold, $X_3$ is $1/2-$ anticanonically embedded Fano and in general $X_i$ is a $1/(i-1)-$ anticanonically embedded Fano $(i+2)-$ fold. Since such polarized varieties are smooth points of their respective Hilbert schemes, we have a chain  $\mathcal{H}_{k+1} \to \mathcal{H}_{k} \to ... \to \mathcal{H}_1 = \mathcal{H}_Y^1$ of non-empty irreducible components of the Hilbert schemes containing the extensions where the maps are given by taking hyperplane sections that send $X_{k+1}$ to $X_{k}$. We show that the map $\mathcal{H}_i \to \mathcal{H}_{i-1}$ is surjective for $2 \leq i \leq k+1$. Denote by $G = Y_{k+1} \supset Y_{k} \supset ... \supset Y_1 = Y$, the chain of extensions of $Y$ and recall that $X_i$ is obtained as a general member of a one parameter family of smoothings of ribbons $\widetilde{Y}_i \in \mathcal{H}_i$ on $Y_i$ with conormal bundle $\mathcal{O}_{Y_i}(-1)$. By standard deformation theory arguments (for example see \cite[Section $5.4$]{CLM93}), it is enough to show that the map $H^1(T_{X_i}) \to H^1(T_{X_i}|_{X_{i-1}})$ is surjective. Hence it is enough to show that $H^2(T_{X_i}(-1)) = 0$. Now $H^2(T_{X_i}(-1)) = H^i(\Omega_{X_i}(K_{X_i} \otimes \mathcal{O}_{X_i}(1)) = 0$ for $i \geq 3$ by Nakano vanishing theorem since $K_{X_i} \otimes \mathcal{O}_{X_i}(1) = (i-2/i-1)K_{X_i}$ and is hence negative ample for $i \geq 3$. For $i = 2$, we need to show that $H^2(\Omega_{X_2}) = 0$. Now $X_2$ is obtained as a smoothing of a ribbon over $Y_2$ with conormal bundle $\mathcal{O}_{Y_2}(-1)$ since the composition of the map $\varphi$ as below \\

\begin{center}
\begin{tikzcd}
    W_2 \arrow[dr, "\varphi"] \arrow[d, "\pi"] & \\
    Y_2 \arrow[r, hook, "|\mathcal{O}_{Y_2}(1)|"] & \mathbb{P}^{M+1}
\end{tikzcd}
\end{center}

deforms to an embedding, where $\pi$ is a double cover branched along a smooth member of $|\mathcal{O}_{Y_2}(2)|$. \textcolor{black}{Note such a smooth member exists by Bertini's theorem as $|\mathcal{O}_{Y_2}(2)|$ is very ample.} Hence $X_2$ is a deformation of $W_2$ and hence it is enough to show that $H^2(\Omega_{W_2}) = 0$. By \cite{Par91}, we need to show that $H^2(\Omega_{Y_2}) = H^2(\Omega_{Y_2}(\textrm{log}(B) \otimes \mathcal{O}_{Y_2}(-1)) = 0$. Since $Y_2$ is a Fano fourfold which is a linear section of a Mukai variety, we have by Lefschetz theorem on Hodge numbers, that $h^{1,2}(Y_2) = 0$. The latter vanishing amounts to showing $H^2(\Omega_{Y_2} \otimes \mathcal{O}_{Y_2}(-1) = H^2(\mathcal{O}_{Y_2}(-1)|_B) = 0$. The former is once again zero due to Nakano vanishing theorem. For the latter we need to show that $H^2(\mathcal{O}_{Y_2}(-1)) = H^3(\mathcal{O}_{Y_2}(-1) \otimes \mathcal{O}_{Y_2}(-B)) = 0$, both of which follows from Kodaira vanishing theorem. \QEDA  


\begin{remark}\label{non-extendability of 1.12-1.15}
 For the rest of the families, i.e, $1.10, 1.12, 1.13-1.15$, the anticanonical embedding of the corresponding Fano-varieties are not smoothly extendable. But some of them are extendable to singular arithmetically Gorenstein normal varieties (see \cite{CLM98}, \cite{CDS20}). To apply the above method, one therefore must have a theory of existence and smoothing of embedded ribbons on singular arithmetically Gorenstein normal varieties. 
\end{remark}

\subsection{Non-extendability of general elements of $\mathcal{H}^Y_l$ for higher values of $l$ and their canonical surface sections}
In this section we study the non-extendability of the general members of $\mathcal{H}^Y_l$ for higher values of $l$. We state and prove the following theorem, which is our second main result on the extendability of general Calabi-Yau threefolds of $\mathcal{H}_l^Y$ and their canonical surface sections.


\begin{theorem}\label{Picard rank one}
\textcolor{black}{Let the notation be as in Proposition \ref{existence and smoothing of ribbons}}
    \begin{itemize}
        \item[(1)]  Then the extendability for a general smooth Calabi-Yau threefold of the Hilbert component $\mathcal{H}^Y_l$ is listed in the fourth column of table \hyperref[t02]{2}. In particular for each of the $17$ different deformation types of $Y$, we give an effective value $l_Y$, such that for $l \geq l_Y$, a smooth Calabi-Yau threefold $X_t \subset \mathbb{P}^{N_l}$ for a general $t$ in $\mathcal{H}^Y_l$ has $\alpha(X_t) = 0$ and is hence not extendable. Consequently, the cone over $X_t$ in $\mathbb{P}^{N_l+1}$ is not smoothable in $\mathbb{P}^{N_l+1}$.
        \item[(2)] For any $Y$ and $l \geq l_Y$, the cohomology groups $H^1(N_{X_t/\mathbb{P}^{N_l}})(-k) = 0$ for $k \geq l$ and $t \in \mathcal{H}^Y_l$ general.
        \item[(3)] Consequently for each deformation type $Y$ and $l \geq l_Y$, the general canonical surface sections $S_t \subset \mathbb{P}^{N_l-1}$ of \textcolor{black}{the} general Calabi-Yau threefolds $X_t \subset \mathbb{P}^{N_l}$ of $\mathcal{H}^Y_l$ have $\alpha(S_t) = 1$, is exactly $1-$ extendable and form a unique irreducible component of the Hilbert scheme $\mathcal{S}^Y_l$ with a dominant map $\mathcal{H}^Y_l \to \mathcal{S}^Y_l$. 
        \item[(4)] For $l \geq l_Y$, $\mathcal{H}^Y_l$ is the only irreducible component of the Hilbert scheme parameterizing linearly normal Calabi-Yau threefolds that \textcolor{black}{dominates} the component $\mathcal{S}^Y_l$ by taking hyperplane sections. For a general canonical surface $S_t \subset \mathbb{P}^{N_l-1} \in \mathcal{S}^Y_l$, the Calabi-Yau threefolds in $\mathbb{P}^{N_l}$ containing $S_t$ as a hyperplane section form an irreducible family. 
    \end{itemize}
\end{theorem}

\begin{remark}\label{KLM11}
Applying \cite[Remark $3.1$]{KLM11}, where the authors apply Zak-L\'vovsky's theorem for veronese embeddings induced by multiples of a very ample line bundle, one can deduce the following weaker bounds for $l_Y$ as follows: 
\begin{itemize}
    \item[(1)] $1.1, l_Y \geq 15$
    \item[(2)] $1.2, l_Y \geq 5$
    \item[(3)] $1.3-1.10, l_Y \geq 5$
    \item[(4)] $1.11, l_Y \geq 15$
    \item[(5)] $1.12-1.15, l_Y \geq 10$
    \item[(6)] $1.16, l_Y \geq 15$
    \item[(7)] $1.17, l_Y \geq 20$
\end{itemize}
In the following table, recall that by Proposition \ref{existence and smoothing of ribbons} and Remark \ref{values of j}, $j$ is the least positive integer such that $lA_t$ is very ample. It is to be noted that apart from the stronger bounds, we also show the canonical surface sections are $1-$ extendable and the uniqueness of the component of the Hilbert scheme of Calabi-Yau threefolds containing the surfaces as their hyperplane sections. We also remark that the gap between the starting value of $l$ in the table and $j+1$ is due to the fact that the upper bound $\beta$ for $\alpha(X_t)$ in Theorem \ref{ZL for ribbons} satisfies $\beta \geq N_l$ and consequently the bound is inconclusive to determine the non-extendability of $X_t \subset \mathbb{P}^{N_l}$. 
\end{remark}

\pagebreak

\begin{center}\phantomsection \label{t02}
\resizebox{7cm}{!}{\begin{tabular}{c|c|c|c|c|c}
\hline
\makecell{Deformation type \\ of Fano threefold \\ $Y$} & \makecell{\textcolor{black}{Value of} \\ \textcolor{black}{$j$}} & \makecell{Value of \\ $l$} & \makecell{extendability of \\ a general Calabi-Yau \\ threefold $X_t$ \\ in $\mathcal{H}_l^Y$}  & \makecell{extendability of \\ a smooth surface \\ section $S_t$ of $X_t$}  & $H^1(N_{X_t/\mathbb{P}^{N_l}}(-l))$ \\
\hline\hline
$1.1$ & \textcolor{black}{$3$} & $ l = 6$  & \makecell{$\alpha(X_t) \leq 1$ \\ hence not $2$-extendable}  & & \\
 \hline
$1.1$ & \textcolor{black}{$3$} & $l \geq 7$  &  \makecell{$\alpha(X_t) = 0$ \\ hence not extendable}  & \makecell{$\alpha(S_t) = 1$ \\ and hence extendable \\ but not $2$-extendable}  & $0$ \\
 \hline
$1.2$ & \textcolor{black}{$1$} & $ l = 4$  & \makecell{$\alpha(X_t) \leq 1$ \\ hence not $2$-extendable}  & & \\
 \hline
$1.2$ & \textcolor{black}{$1$} & $l \geq 5$  &  \makecell{$\alpha(X_t) = 0$ \\ hence not extendable}  & \makecell{$\alpha(S_t) = 1$ \\ and hence extendable \\ but not $2$-extendable}  & $0$\\
 \hline
$1.3$ & \textcolor{black}{$1$} & $ l = 3$  & \makecell{$\alpha(X_t) \leq 1$ \\ hence not $2$-extendable}  & & \\
 \hline
$1.3$ & \textcolor{black}{$1$} &  $l \geq 4$   & \makecell{$\alpha(X_t) = 0$ \\ hence not extendable} & \makecell{$\alpha(S_t) = 1$ \\ and hence extendable \\ but not $2$-extendable}  & $0$\\
 \hline 
$1.4$ &  \textcolor{black}{$1$} & $l = 2$   & \makecell{$\alpha(X_t) \leq 4$ \\ hence not $5$-extendable} & &  \\
\hline 
$1.4$ & \textcolor{black}{$1$} & $l \geq 3$   & \makecell{$\alpha(X_t) = 0$ \\ hence not extendable} & \makecell{$\alpha(S_t) = 1$ \\ and hence extendable \\ but not $2$-extendable}  & $0$\\
\hline
$1.5$ & \textcolor{black}{$1$} & $l = 2$   & \makecell{$\alpha(X_t) \leq 2$ \\ hence not $3$-extendable} & &  \\
 \hline
 $1.5$ & \textcolor{black}{$1$} & $l \geq 3$   & \makecell{$\alpha(X_t) = 0$ \\ hence not extendable} & \makecell{$\alpha(S_t) = 1$ \\ and hence extendable \\ but not $2$-extendable}  & $0$\\
 \hline
$1.6-1.10$ & \textcolor{black}{$1$} & $l = 2$   & \makecell{$\alpha(X_t) \leq 1$ \\ hence not $2$-extendable} & & \\
 \hline
 $1.6-1.10$ & \textcolor{black}{$1$} & $l \geq 3$   & \makecell{$\alpha(X_t) = 0$ \\ hence not extendable} & \makecell{$\alpha(S_t) = 1$ \\ and hence extendable \\ but not $2$-extendable}  & $0$\\
 \hline
$1.11$ & \textcolor{black}{$3$} & $l = 4$  & \makecell{$\alpha(X_t) \leq 12$ \\ hence not $13$-extendable} & & \\
 \hline
 $1.11$ & \textcolor{black}{$3$} & $l = 5$  & \makecell{$\alpha(X_t) \leq 3$ \\ hence not $4$-extendable} & & \\
 \hline
 $1.11$ & \textcolor{black}{$3$} & $l = 6$  & \makecell{$\alpha(X_t) \leq 1$ \\ hence not $2$-extendable} & & \\
 \hline
$1.11$ & \textcolor{black}{$3$}  & $l \geq 7$  & \makecell{$\alpha(X_t) = 0$ \\ hence not extendable} & \makecell{$\alpha(S_t) = 1$ \\ and hence extendable \\ but not $2$-extendable}  & $0$ \\
 \hline
$1.12$ & \textcolor{black}{$2$} & $l = 3$  & \makecell{$\alpha(X_t) \leq 8$ \\ hence not $9$-extendable} & & \\
 \hline
$1.12$ & \textcolor{black}{$2$} & $l = 4$  & \makecell{$\alpha(X_t) \leq 2$ \\ hence not $3$-extendable} & & \\
 \hline
$1.12$ & \textcolor{black}{$2$} & $l \geq 5$  & \makecell{$\alpha(X_t) = 0$ \\ hence not extendable} & \makecell{$\alpha(S_t) = 1$ \\ and hence extendable \\ but not $2$-extendable}  & $0$ \\
 \hline
$1.13-1.14$ & \textcolor{black}{$2$} &  $l = 3$  & \makecell{$\alpha(X_t) \leq 6$ \\ hence not $7$-extendable} & & \\
 \hline
 $1.13-1.14$ & \textcolor{black}{$2$} & $l = 4$  & \makecell{$\alpha(X_t) \leq 1$ \\ hence not $2$-extendable} & & \\
 \hline
$1.13-1.14$ & \textcolor{black}{$2$} &  $l \geq 5$  & \makecell{$\alpha(X_t) = 0$ \\ hence not extendable} & \makecell{$\alpha(S_t) = 1$ \\ and hence extendable \\ but not $2$-extendable}  & $0$ \\
 \hline
$1.15$ & \textcolor{black}{$2$} &  $l = 3$  & \makecell{$\alpha(X_t) \leq 6$ \\ hence not $7$-extendable} & & \\
\hline
$1.15$ & \textcolor{black}{$2$} &  $l = 4$  & \makecell{$\alpha(X_t) \leq 1$ \\ hence not $2$-extendable} & & \\
\hline
$1.15$ & \textcolor{black}{$2$} &  $l \geq 5$  & \makecell{$\alpha(X_t) = 0$ \\ hence not extendable} & \makecell{$\alpha(S_t) = 1$ \\ and hence extendable \\ but not $2$-extendable}  & $0$ \\
\hline
$1.16$ & \textcolor{black}{$3$} &  $l = 4$  & \makecell{$\alpha(X_t) \leq 20$ \\ hence not $21$-extendable} & & \\
\hline
$1.16$ & \textcolor{black}{$3$} & $l = 5$  & \makecell{$\alpha(X_t) \leq 5$ \\ hence not $6$-extendable} & & \\
\hline
$1.16$ & \textcolor{black}{$3$} & $l \geq 6$  & \makecell{$\alpha(X_t) = 0$ \\ hence not extendable} & \makecell{$\alpha(S_t) = 1$ \\ and hence extendable \\ but not $2$-extendable}  & $0$ \\
\hline
$1.17$ & \textcolor{black}{$4$} & $l = 5$  & \makecell{$\alpha(X_t) \leq 35$ \\ hence not $36$-extendable} & & \\
\hline
$1.17$ & \textcolor{black}{$4$} & $l = 6$  & \makecell{$\alpha(X_t) \leq 15$ \\ hence not $16$-extendable} & & \\
\hline
$1.17$ & $4$ & $l = 7$  & \makecell{$\alpha(X_t) \leq 4$ \\ hence not $5$-extendable} & & \\
\hline
$1.17$ & $4$ & $l = 8$  & \makecell{$\alpha(X_t) \leq 1$ \\ hence not $2$-extendable} & & \\
\hline
$1.17$ & $4$  & $l \geq 9$  & \makecell{$\alpha(X_t) = 0$ \\ hence not extendable} & \makecell{$\alpha(S_t) = 1$ \\ and hence extendable \\ but not $2$-extendable}  & $0$\\
\hline

\end{tabular}}
\captionof{table}{}
\end{center}

\medskip


\medskip

\noindent\textit{Proof of Theorem \ref{Picard rank one}.}  
\textcolor{black}{We briefly indicate the strategy of how we prove $(1)-(4)$ in the statement of the theorem and show that it eventually boils down to computing upper bounds for $h^i(T_Y(-k))$ for $k \geq l$ and $i = 1,2$:}
\begin{itemize}
\item[(1)] To prove statement $(1)$, we need to compute $\alpha(X_t)$ for $t \in \mathcal{H}^Y_l$ general. By Proposition \ref{existence and smoothing of ribbons}, $(4)$, we calculate $\alpha(\widetilde{Y})$ for the Calabi-Yau ribbons $\widetilde{Y} \subset \mathbb{P}^{N_l}$ embedded by $\mathcal{O}_{\widetilde{Y}}(1)$ where $\mathcal{O}_{\widetilde{Y}}(1)|_Y = lH$. So we apply Theorem \ref{ZL for ribbons} with $\mathcal{O}_{\widetilde{Y}}(1) = \widetilde{H}$. Note that since $\widetilde{Y}$ is a Calabi-Yau ribbon, we have that $L = K_Y = -jH$. Theorem \ref{ZL for ribbons} is effective when $l > j$ where $j$ is the index of the Fano threefold.  Note that if 
\begin{enumerate}\label{calculate alpha}
    \item[(a)] $j < l < 2j$, we have $\beta \leq h^1(T_Y(-(l+j)) + h^1(T_Y(-l)) - h^0(T_Y(-l)) + h^0(\mathcal{O}_Y(-l+2j)$)
    \item[(b)] $l = 2j$, we have $\beta \leq h^1(T_Y(-(l+j)) + h^1(T_Y(-l)) - h^0(T_Y(-l)) + 1$
    \item[(c)] $l > 2j$, we have $\beta \leq h^1(T_Y(-(l+j)) + h^1(T_Y(-l)) - h^0(T_Y(-l))$
\end{enumerate}
So we are down to computing $H^1(T_Y(-l))$ for $l > j$ for the $17$ different deformation classes and for each family we will find an $l_Y > j$ such that the above cohomology group vanishes for $l \geq l_Y$.

\smallskip

\item[(2)] To prove statement $(2)$, in the case $l \geq l_Y$, we need to show the vanishing of $H^1(N_{X_t/\mathbb{P}^{N_l}})(-k) = 0$ for $k \geq l$ and $t \in \mathcal{H}^Y_l$ general. We prove the vanishing of $H^1(N_{\widetilde{Y}/\mathbb{P}^{N_l}})(-k\widetilde{H})$ with $\widetilde{H}|_Y = lH$, \textcolor{black}{as before}, \textcolor{black}{and apply} Proposition \ref{cohomology of twists of normal bundles}. Now note that all the conditions of part $(1)$ of Proposition \ref{cohomology of twists of normal bundles} are satisfied since intermediate cohomologies of any ample line bundle or the structure sheaf of a Fano variety vanishes. So we need to show that $H^1(T_Y(-(kl+j))) = H^2(T_Y(-kl)) = 0$ for $l \geq l_Y$ and $1 \leq k \leq 3$. Note that from the proof of part $(1)$ we already know $H^1(T_Y(-(kl+j))) = 0$ since $k \geq 1$ and $j \geq 1$. For the second vanishing we will show that $H^2(T_Y(-l)) = 0$ for $l \geq l_Y$.

\smallskip

\item[(3)] Part $(3)$ now follows from Lemma \ref{hyperplane sections fill out entire moduli}.

\smallskip

\item[(4)] We first show that a general Calabi-Yau in $\mathcal{H}_l^Y$ and hence a general canonical surface section in $\mathcal{S}_l^Y$ for $l \geq l_Y$ is projectively normal. By \cite[Corollary $1.1$]{GP97}, for a Calabi-Yau threefold $X$ and an ample and base point free line bundle $B$ on $X$, $lB$ is projectively normal for $l \geq 4$. For a general $X_t \in \mathcal{H}_l^Y$, the generator $A_t$ of the Picard group is base point free and from the tables \hyperref[t01]{1}, \hyperref[t02]{2}, we see that $l_Y \geq 4$ for all deformation type $Y$ other than $1.4-1.10$. Now once again, by \cite[Theorem $1$]{GP97}, for a Calabi-Yau threefold $X$ and an ample and base point free line bundle $B$ on $X$, $lB$ is projectively normal for $l \geq 3$, unless the morphism by $|B|$ maps $X$, $2:1$ onto $\mathbb{P}^3$. In our case therefore for families $1.4-1.10$, $h^0(lA_t) > 4$ by Proposition \ref{existence and smoothing of ribbons}, part $(2)$ and hence are projectively normal for $l \geq 3$. By \cite[remark $7.6$ (iii)]{Pin74}, the cone $C(Z) \subset \mathbb{P}^{N_l}$ over a general hyperplane section $Z \subset \mathbb{P}^{N_l-1}$ lies in every irreducible component of the Hilbert scheme in $\mathbb{P}^N$ that dominates $\mathcal{S}_l^Y$. But now by Lemma \ref{smoothness of Hilbert point of the cone} the Hilbert point of such a cone is smooth. Hence part $(4)$ now follows.

\end{itemize}

\smallskip

For the rest of the proof, \textcolor{black}{as noted towards the end of (1) and (2) above, we need to} estimate $h^i(T_Y(-k))$ for $i = 1,2$ for each of the $17$ deformation types of Fano-threefolds $Y$ and compute $\alpha(X)$ according to formulas (a), (b), (c) \textcolor{black}{in} $(1)$ above. \textcolor{black}{To compute the cohomology, we use the Borel-Weil-Bott theorem (see for example \cite[Section $2.6$]{Kuz05}) and some established vanishing theorems on the cohomology of twisted holomorphic forms on Grassmannians in \cite{Sno86} (see also \cite[Theorem $3.6$]{FM21}). We further use Pieri's rule as stated in \cite{Pieri1893} or \cite[Theorem $2.3$]{AQ18}.} \textcolor{black}{In the computations to follow, we assume $ i = 1,2$.}
\smallskip

\textcolor{black}{\bf{We prove for $Y$ in family $1.7$.}} Let $G_n$ denote the Grassmannian Gr$(2,n+2)$ so that we have Gr $(2,6)$ is denoted by $G_4$. $G_4$ is a Fano variety of dimension $8$ and index $6$ (see for ex \cite{Man15}, Section $(2)$) while its Picard group is generated by the line bundle $\mathcal{O}_{G_4}(1)$ corresponding to the Plucker embedding. Hence $K_{G_4} = \mathcal{O}_{G_4}(-6)$. On the other hand $Y$ is the intersection of $G_4$ in its Plucker embedding by a codimension $5$ subspace. Therefore, $K_Y = \mathcal{O}_Y(-1)$ where $\mathcal{O}_Y(1)$ is the restriction of $\mathcal{O}_{G_4}(1)$ to $Y$. Let us show $H^i(T_Y(-l)) = 0$ for $l \geq 2$. Noting that the normal bundle of $Y$ inside $G_4$ is $\mathcal{O}_Y(1)^{\oplus 5}$ we have the exact sequence 
$$0 \to T_Y \to T_{G_4}|_Y \to \mathcal{O}_Y(1)^{\oplus 5} \to 0$$
Tensoring by $\mathcal{O}_Y(-l)$ we have for $i = 1,2$ 
$$H^{i-1}(\mathcal{O}_Y(1-l))^{\oplus 5} \to H^i(T_Y(-l)) \to H^i(T_{G_4}|_Y(-l))$$
For $l \geq 2$, we have that $H^{i-1}(\mathcal{O}_Y(1-l)) = 0$. Now let $\mathcal{I}_Y$ be the ideal sheaf of $Y$ inside $G_4$. Then we have that 
$$0 \to \mathcal{I}_Y \to \mathcal{O}_{G_4} \to \mathcal{O}_Y \to 0$$
Tensoring by $T_{G_4}(-l)$ and taking cohomology we have
$$H^i(T_{G_4}(-l)) \to H^i(T_{G_4}|_Y(-l)) \to H^{i+1}(T_{G_4}(-l) \otimes \mathcal{I}_Y)$$
Recall that by Borel-Bott-Weil theorem, $H^q(\Omega_{G_4}(t)) = 0$ for $(q,t) \neq (1,0)$ and $q \geq 1$. Then $H^i(T_{G_4}(-l)) = H^{8-i}(\Omega_{G_4}(l-6)) = 0$ for $i = 1,2$. Now if $W$ is the five dimensional vector space generated by $s_0,..s_4$ , the sections of $\mathcal{O}_{G_4}(1)$ that cuts out $Y$, then the ideal sheaf $\mathcal{I}_Y$ has a Koszul resolution 
$$0 \to \bigwedge^5 W \otimes \mathcal{O}_{G_4}(-5)\to...\to \bigwedge^2 W \otimes \mathcal{O}_{G_4}(-2) \to W \otimes \mathcal{O}_{G_4}(-1) \to \mathcal{I}_Y \to 0$$

So if we show $H^{i+k+1}(T_{G_4}(-l-k-1))$ for $i = 1,2$ and $k = 0,...,4$, then we are done. This amounts to showing the vanishing of $H^{7-i-k}(\Omega_{G_4}(l+k-5)) = 0$ for $i = 1,2$ and $k = 0,...,4$. But this once again follows from Borel-Bott-Weil theorem. Therefore our result now follows. \par

\smallskip

\textcolor{black}{\bf{The proof for family $1.5$}} follows along the same lines as $1.7$. In this case, $Y$ is cut out inside the plucker embedding of $G_3$ by two hyperplanes and a quadric. We have $K_{G_3} = \mathcal{O}_{G_3}(-5)$ and $K_Y = \mathcal{O}_{Y}(-1)$. Noting that the normal bundle of $Y$ inside $G_3$ is $\mathcal{O}_Y(1)^{\oplus 2} \bigoplus \mathcal{O}_Y(2)$ we have the exact sequence 
$$0 \to T_Y \to T_{G_3}|_Y \to \mathcal{O}_Y(1)^{\oplus 2} \bigoplus \mathcal{O}_Y(2) \to 0$$
Tensoring by $\mathcal{O}_Y(-l)$ and taking cohomology we have for $i = 1,2$
$$H^{i-1}(\mathcal{O}_Y(1-l))^{\oplus 2} \bigoplus H^{i-1}(\mathcal{O}_Y(2-l)) \to H^i(T_Y(-l)) \to H^i(T_{G_3}|_Y(-l))$$
For $l \geq 3$, we have that $H^{i-1}(\mathcal{O}_Y(1-l)) = H^{i-1}(\mathcal{O}_Y(2-l)) = 0$ while for $l = 2$, $H^0(\mathcal{O}_Y(1-l)) = 0$ and $H^0(\mathcal{O}_Y(2-l)) = 1$ . Now let $\mathcal{I}_Y$ be the ideal sheaf of $Y$ inside $G_3$. Then we have that 
$$0 \to \mathcal{I}_Y \to \mathcal{O}_{G_3} \to \mathcal{O}_Y \to 0$$
Tensoring by $T_{G_3}(-l)$ and taking cohomology we have
$$H^i(T_{G_3}(-l)) \to H^i(T_{G_3}|_Y(-l)) \to H^{i+1}(T_{G_3}(-l) \otimes \mathcal{I}_Y)$$
Then $H^i(T_{G_3}(-l)) = H^{5-i}(\Omega_{G_3}(l-5)) = 0$ by Borel-Bott-Weil theorem. Now the ideal sheaf $\mathcal{I}_Y$ has a resolution 
$$0 \to \mathcal{O}_{G_3}(-4) \to \mathcal{O}_{G_3}(-2)^{\oplus 2} \bigoplus \mathcal{O}_{G_3}(-3) \to \mathcal{O}_{G_3}(-1)^{\oplus 2} \bigoplus \mathcal{O}_{G_3}(-2) \to \mathcal{I}_Y \to 0$$

Then it is enough to show show $H^{i+1}(T_{G_3}(-1-l)) = H^{i+1}(T_{G_3}(-2-l)) = H^{i+2}(T_{G_3}(-2-l)) = H^{i+2}(T_{G_3}(-3-l)) = H^{i+3}(T_{G_3}(-4-l)) = 0$ for $l \geq 1$ and this follows by Borel-Bott-Weil theorem. So we have that $h^i(T_Y(lK_Y)) = 0$ for $l \geq 3$, $i = 1,2$ and \textcolor{black}{$h^1(T_Y(lK_Y)) \leq 1$} for $l = 2$. Hence our result follows. \par

\smallskip

\textcolor{black}{\bf{The proof for family $1.15$}} follows along the same lines as $1.5$. In this case, $Y$ is cut out inside the Plucker embedding of $G_3$ by three hyperplanes. We have $K_{G_3} = \mathcal{O}_{G_3}(-5)$ and $K_Y = \mathcal{O}_{Y}(-2)$. Noting that the normal bundle of $Y$ inside $G_3$ is $\mathcal{O}_Y(1)^{\oplus 3} $ we have the exact sequence 
$$0 \to T_Y \to T_{G_3}|_Y \to \mathcal{O}_Y(1)^{\oplus 3}  \to 0$$
Tensoring by $\mathcal{O}_Y(-l)$ with $l \geq 2$ and taking cohomology we have for $i = 1,2$, 
$$H^{i-1}(\mathcal{O}_Y(1-l))^{\oplus 3}  \to H^i(T_Y(-l)) \to H^i(T_{G_3}|_Y(-l))$$
For $l \geq 2$, we have that $H^{i-1}(\mathcal{O}_Y(1-l)) = 0$. Now let $\mathcal{I}_Y$ be the ideal sheaf of $Y$ inside $G_3$. Then we have that 
$$0 \to \mathcal{I}_Y \to \mathcal{O}_{G_3} \to \mathcal{O}_Y \to 0$$
Tensoring by $T_{G_3}(-2l)$ and taking cohomology we have
$$H^i(T_{G_3}(-l)) \to H^i(T_{G_3}|_Y(-l)) \to H^{i+1}(T_{G_3}(-l) \otimes \mathcal{I}_Y)$$
Then $H^i(T_{G_3}(-l)) = H^{6-i}(\Omega_{G_3}(l-5)) = 0$ by Borel-Bott-Weil theorem. Now the ideal sheaf $\mathcal{I}_Y$ has a resolution 
$$0 \to \mathcal{O}_{G_3}(-3) \to \mathcal{O}_{G_3}(-2)^{\oplus 3} \to \mathcal{O}_{G_3}(-1)^{\oplus 3} \to \mathcal{I}_Y \to 0$$

As before, it is enough to show $H^{i+1}(T_{G_3}(-1-l)) = H^{i+2}(T_{G_3}(-2-l)) = H^{i+3}(T_{G_3}(-3-l)) = 0$ and this follows for any $l$ by Borel-Bott-Weil theorem. So we have that $h^i(T_Y(-l)) = 0$ for $l \geq 2$ and $i = 1,2$ and hence our result follows.

\textcolor{black}{}{\textcolor{black}{\bf{The proof for family $1.6$}}. In this case we use the description that $Y$ is the zero locus of a section of the vector bundle $E := \mathcal{U}^*(1) \oplus \mathcal{O}_{G_3}(1)$ where $\mathcal{U}$ is the tautological rank two vector bundle of $G_3$. Recall that $K_{G_3} = \mathcal{O}_{G_3}(-5)$ and $N_{Y/G_{3}} = \mathcal{U}^*_Y(1) \oplus \mathcal{O}_{Y}(1) = E_Y$ where $\mathcal{U}_Y$ is the restriction of $\mathcal{U}$ to $Y$. Therefore, by adjunction we have that $K_Y = K_{G_3} \otimes \textrm{det}(\mathcal{U}^*_Y)(3)$. Recall that $\textrm{det}(\mathcal{U}^*) = \mathcal{O}_{G_3}(1)$. 
Hence $K_Y = \mathcal{O}_Y(-1)$. Twisting the sequence 
$$0 \to T_Y \to T_{G_3}|_Y \to \mathcal{U}^*_Y(1) \oplus \mathcal{O}_{Y}(1) \to 0$$ 
by $\mathcal{O}_Y(-l)$ and taking cohomology, we have for $i = 1,2$ 
$$H^{i-1}(\mathcal{U}^*_Y(1-l) \oplus H^{i-1}(\mathcal{O}_{Y}(1-l))) \to H^i(T_Y(-l)) \to H^{i}(T_{G_3}|_Y(-l))$$
We first show the vanishing of $H^{i-1}(\mathcal{U}^*_Y(1-l) \oplus H^{i-1}(\mathcal{O}_{Y}(1-l)))$. Clearly, $H^{i-1}(\mathcal{O}_{Y}(1-l))) = 0$ for $l \geq 2$ and $i = 1,2$. Now note that there exist an exact sequence
$$0 \to \mathcal{U} \to \mathcal{O}_{G_3}^{\oplus 5} \to \mathcal{Q} \to 0$$ where $\mathcal{Q}$ is the tautological quotient bundle. Restricting to $Y$, dualizing, tensoring by $\mathcal{O}_Y(1-l)$ and taking cohomology we have that 
$$H^{i-1}(\mathcal{O}_Y(1-l))^{\oplus 5} \to H^{i-1}(\mathcal{U}^*_Y(1-l)) \to H^i(\mathcal{Q}^*_Y(1-l))$$

Clearly for $l \geq 2$, we have that $H^{i-1}(\mathcal{O}_Y(1-l) = 0$. Tensoring the defining exact sequence of $\mathcal{O}_Y$ by $\mathcal{Q}^*(1-l)$ and taking cohomology we have that 
$$H^{i}(\mathcal{Q}^*(1-l)) \to H^{i}(\mathcal{Q}_Y^*(1-l)) \to H^{i+1}(\mathcal{Q}^* \otimes \mathcal{I}_Y(1-l))$$
We have that $H^i(\mathcal{Q}^*(1-l)) = H^{6-i}(\mathcal{Q}(l-6)) = H^{6-i}(\Sigma^{(l-6, l-6)} \mathcal{U}^* \otimes \Sigma^{(1,0,0)} \mathcal{Q}) = \textcolor{black}{H^{6-i}(\Sigma^{(l-6, l-6)} \mathcal{U}^* \otimes \Sigma^{(0,0,-1)} \mathcal{Q}^*)} = H^{6-i}(L_{(l-6, l-6, 0, 0, -1)})$.
Setting $\alpha = \textcolor{black}{(l-6, l-6, 0, 0, -1)}$ and $\rho = (5,4,3,2,1)$, we have that $\alpha+\rho  = \textcolor{black}{(l-1, l-2, 3, 2, 0)}$. For any $2 \leq l \leq  5$, the entries of $\alpha+\rho$ are not distinct while if $l \geq  6$, the entries are in strictly decreasing order and consequently the permutation that arranges the sequence in strictly decreasing order is identity whose length is zero. Hence $ H^{6-i}(L_{(l-6, l-6, 0, 0, -1)}) = 0$ by the Borel-Weil-Bott theorem (see \cite[Section $2.6$]{Kuz05}).
Now consider the Koszul resolution of the ideal sheaf $\mathcal{I}_Y$
$$ 0 \to \textrm{det}(\mathcal{U})(-3) \to \bigwedge^2 E^* \to \mathcal{U}(-1) \oplus \mathcal{O}_{G_3}(-1) \to \mathcal{I}_Y \to 0 $$
Tensoring with $\mathcal{Q}^*(1-l)$, we have 
$$ 0 \to \mathcal{Q}^*(-3-l) \to \mathcal{Q}^*(1-l) \otimes \bigwedge^2 E^* \to \mathcal{Q}^* \otimes \mathcal{U}(-l) \oplus \mathcal{Q}^* (-l) \to \mathcal{Q}^* \otimes \mathcal{I}_Y (1-l) \to 0 $$
We have that $H^{i+1}(\mathcal{Q}^* (-l)) = H^{5-i}(\mathcal{Q}(l-5)) = H^{5-i}(\Sigma^{(l-5,l-5)} \mathcal{U}^* \otimes \Sigma^{(1,0,0)}\mathcal{Q}) = H^{5-i}(\Sigma^{(l-5,l-5)} \mathcal{U}^* \otimes \Sigma^{(0,0,-1)}\mathcal{Q}^*) = H^{5-i}(L_{(l-5,l-5,0,0,-1)}) = 0$ for $l \geq 2$ by the Borel-Bott-Weil theorem as above. 

$H^{i+1}(\mathcal{Q}^* \otimes \mathcal{U}(-l)) = H^{5-i}(\mathcal{U}^*(l-5) \otimes \mathcal{Q}) = H^{5-i}((\Sigma^{(l-4,l-5)} \mathcal{U}^* \otimes \Sigma^{(1,0,0)}\mathcal{Q})) = H^{5-i}((\Sigma^{(l-4,l-5)} \mathcal{U}^* \otimes \Sigma^{(0,0,-1)}\mathcal{Q}^*)) = H^{5-i}(L_{(l-4,l-5,0,0,-1)}) = 0$ for $l \geq 2$ by the Borel-Bott-Weil theorem. 

$H^{i+3}(\mathcal{Q}^*(-3-l)) = H^{3-i}(\mathcal{Q}(l-2)) = H^{3-i}((\Sigma^{(l-2,l-2)} \mathcal{U}^* \otimes \Sigma^{(1,0,0)}\mathcal{Q})) = H^{3-i}((\Sigma^{(l-2,l-2)} \mathcal{U}^* \otimes \Sigma^{(0,0,-1)}\mathcal{Q}^*)) = H^{3-i}(L_{(l-2,l-2,0,0,-1)}) = 0$ for $l \geq 2$ by the Borel-Bott-Weil theorem.

Now consider $H^{i+2}(\mathcal{Q}^*(1-l) \otimes \bigwedge^2 E^*)$. Recall that $E^* = \mathcal{U}(-1) \oplus \mathcal{O}_{G_3}(-1)$. We have 
$$0 \to \mathcal{U}(-1) \to E^* \to \mathcal{O}_{G_3}(-1) \to 0$$ and hence
$$ 0 \to \bigwedge^2 (\mathcal{U}(-1)) \to \bigwedge^2 E^*\to \mathcal{U}(-2) \to 0 $$
Now $\bigwedge^2 (\mathcal{U}(-1)) = \textrm{det}(\mathcal{U})(-2) = \mathcal{O}_{G_3}(-3)$. So tensoring the above sequence by $\mathcal{Q}^*(1-l)$ and taking cohomology, we have
$$H^{i+2}(\mathcal{Q}^*(-2-l)) \to H^{i+2}(\mathcal{Q}^*(1-l) \otimes \bigwedge^2 E^*) \to H^{i+2}(\mathcal{Q}^* \otimes \mathcal{U}(-1-l))$$

We have $H^{i+2}(\mathcal{Q}^*(-2-l)) = H^{4-i}(\mathcal{Q}(l-3)) = H^{4-i}((\Sigma^{(l-3,l-3)} \mathcal{U}^* \otimes \Sigma^{(1,0,0)}\mathcal{Q})) = H^{4-i}((\Sigma^{(l-3,l-3)} \mathcal{U}^* \otimes \Sigma^{(0,0,-1)}\mathcal{Q}^*)) = H^{4-i}(L_{(l-3,l-3,0,0,-1)}) = 0$ for $l \geq 2$ by the Borel-Bott-Weil theorem.

Consider $H^{i+2}(\mathcal{Q}^* \otimes \mathcal{U}(-1-l)) = H^{4-i}(\mathcal{U}^*(l-4) \otimes \mathcal{Q}) = H^{4-i}((\Sigma^{(l-3,l-4)} \mathcal{U}^* \otimes \Sigma^{(1,0,0)}\mathcal{Q})) = H^{4-i}((\Sigma^{(l-3,l-4)} \mathcal{U}^* \otimes \Sigma^{(0,0,-1)}\mathcal{Q}^*)) = H^{4-i}(L_{(l-3,l-4,0,0,-1)})= 0$ for $l \geq 2$ by the Borel-Bott-Weil theorem.. 
This proves that $H^{i+2}(\mathcal{Q}^*(1-l) \otimes \bigwedge^2 E^*) = 0$ for $l \geq 2$ and $i = 1,2$. Therefore $H^{i+1}(\mathcal{Q}^* \otimes \mathcal{I}_Y (1-l)) = 0$ and hence $H^i(\mathcal{Q}^*_Y(1-l)) = 0$. This implies $H^{i-1}(\mathcal{U}_Y^*(1-l)) = 0$ for $l \geq 2$ which in turn gives us the exact sequence 
$$0 \to H^i(T_Y(-l)) \to H^i(T_{G_3}|_Y(-l))$$ for $l \geq 2$. Once again consider the exact sequence 
$$0 \to T_{G_3} \otimes \mathcal{I}_Y (-l) \to T_{G_3}(-l) \to T_{G_3}|_Y(-l) \to 0$$
Taking cohomology, we have
$$H^i(T_{G_3}(-l)) \to H^i(T_{G_3}|_Y(-l)) \to H^{i+1}(T_{G_3} \otimes \mathcal{I}_Y (-l))$$
We have that $H^i(T_{G_3}(-l)) = H^{6-i}(\Omega_{G_3}(l-5)) = 0$. To show that $H^{i+1}(T_{G_3} \otimes \mathcal{I}_Y (-l)) = 0$, we have once again using the Koszul resolution of the ideal sheaf $\mathcal{I}_Y$
$$ 0 \to T_{G_3}(-4-l) \to T_{G_3} \otimes \bigwedge^2 E^*(-l) \to T_{G_3} \otimes \mathcal{U}(-1-l) \oplus T_{G_3}(-1-l) \to T_{G_3} \otimes \mathcal{I}_Y (-l) \to 0 $$

We have $H^{i+1}(T_{G_3}(-1-l)) = H^{i+1}(\mathcal{U}^* \otimes \mathcal{Q} (-1-l)) = H^{5-i}(\mathcal{U} \otimes \mathcal{Q}^* (l-4)) =  H^{5-i}(\Sigma^{(l-4,l-5)} \mathcal{U}^* \otimes \Sigma^{(1,0,0)} \mathcal{Q}^*) = H^{5-i}(L_{(l-4,l-5,1,0,0)}) = 0$ for $l \geq 2$ by the Borel-Bott-Weil theorem. 

\smallskip

Consider $H^{i+1}(T_{G_3} \otimes \mathcal{U}(-1-l))$. Note that $T_{G_3} = \mathcal{U}^* \otimes \mathcal{Q}$ 
and hence $T_{G_3} \otimes \mathcal{U}(-1-l) = \mathcal{U}^* \otimes \mathcal{U} (-1-l) \otimes \mathcal{Q}$. Since $\mathcal{U}$ is a rank two vector bundle we have that $\mathcal{U} = \mathcal{U}^*(-1)$. Hence $H^{i+1}(T_{G_3} \otimes \mathcal{U}(-1-l)) = H^{i+1}(\mathcal{U}^* \otimes \mathcal{U}(-1-l) \otimes \mathcal{Q}) = H^{5-i}(\mathcal{U} \otimes \mathcal{U}^*(l-4) \otimes \mathcal{Q}^*) = H^{5-i}(\mathcal{U}^* \otimes \mathcal{U}^*(l-5) \otimes \mathcal{Q}^*)$. Now by Pieri's rule, $\mathcal{U}^* \otimes \mathcal{U}^* = \Sigma^{(2,0)} \mathcal{U}^* \oplus \Sigma^{(1,1)} \mathcal{U}^*$. Therefore $H^{5-i}(\mathcal{U}^* \otimes \mathcal{U}^*(l-5) \otimes \mathcal{Q}^*) = H^{5-i}(\Sigma^{(l-3,l-5)} \mathcal{U}^* \otimes \Sigma^{(1,0,0)} \mathcal{Q}^*) \oplus H^{5-i}(\Sigma^{(l-4,l-4)} \mathcal{U}^* \otimes \Sigma^{(1,0,0)} \mathcal{Q}^*).$ Both of the terms are zero for $l \geq 2$ by the Borel-Bott-Weil theorem.

\smallskip

Consider $H^{i+2}(T_{G_3} \otimes \bigwedge^2E^*(-l))$. We have $\bigwedge^2 E^*(-l) = \mathcal{O}_{G_3}(-3-l) \oplus \mathcal{U}(-2-l)$. We have $H^{i+2}(T_{G_3}(-3-l)) = H^{4-i}(\Omega_{G_3}(l-2)) = 0$ for $l \geq 2$. The other summand is $H^{i+2}(T_{G_3} \otimes \mathcal{U}(-2-l)) = H^{i+2}(\mathcal{U} \otimes \mathcal{U}^*(-2-l) \otimes \mathcal{Q}) = H^{4-i}(\mathcal{U} \otimes \mathcal{U}^*(l-3) \otimes \mathcal{Q}^*) = H^{4-i}(\mathcal{U}^* \otimes \mathcal{U}^*(l-4) \otimes \mathcal{Q}^*)$. The above has two summands, Similarly $H^{4-i}(\Sigma^{(l-2, l-4)} \mathcal{U}^* \otimes \Sigma^{(1,0,0)} \mathcal{Q}^*) \oplus H^{4-i}(\Sigma^{(l-3, l-3)} \mathcal{U}^* \otimes \Sigma^{(1,0,0)} \mathcal{Q}^*) = 0$ for $l \geq 2$ by the Borel-Bott-Weil theorem. 

\smallskip

Finally , we have $H^{i+3}(T_{G_3}(-4-l)) = H^{i+3}(\mathcal{U}^* \otimes \mathcal{Q} (-4-l)) = H^{3-i}(\mathcal{U} \otimes \mathcal{Q}^* (l-1)) = H^{3-i}(\mathcal{U}^* \otimes \mathcal{Q}^* (l-2)) = H^{3-i}(\Sigma^{(l-1,l-2)} \mathcal{U}^* \otimes \Sigma^{(1,0,0)} \mathcal{Q}^*) = 0$ for $l \geq 2$ 

\smallskip

\textcolor{black}{\bf{We prove for family $1.9$.}} In this case , $Y$ is the zero locus of a section of the vector bundle $E = \mathcal{Q}^*(1) \oplus \mathcal{O}(1)^{\oplus 2}$ on $G_5$. We have that $K_{G_5} = \mathcal{O}_{G_5}(-7)$ while $N_{Y/G_{5}} = \mathcal{Q}_Y^*(1) \oplus \mathcal{O}_Y(1)^{\oplus 2}$. Consequently $K_Y = K_{G_5} \otimes \textrm{det}(\mathcal{Q}_Y)^*(5) \otimes \mathcal{O}_Y(2)$. Considering $\textrm{det}(\mathcal{Q}_Y)^* = \mathcal{O}_Y(-1)$, we have that $K_Y = \mathcal{O}_Y(-1)$. Twisting the sequence 
$$0 \to T_Y \to T_{G_5}|_Y \to \mathcal{Q}_Y^*(1) \oplus \mathcal{O}_Y(1)^{\oplus 2} \to 0$$ 
by $\mathcal{O}_Y(-l)$ and taking cohomology, we have for $i = 1,2$
$$H^{i-1}(\mathcal{Q}_Y^*(1-l) \oplus H^{i-1}(\mathcal{O}_{Y}(1-l))^{\oplus 2}) \to H^i(T_Y(-l)) \to H^i(T_{G_5}|_Y(-l))$$
We first show the vanishing of $H^{i-1}(\mathcal{Q}_Y^*(1-l) \oplus H^{i-1}(\mathcal{O}_{Y}(1-l))^{\oplus 2})$. Clearly, $H^{i-1}(\mathcal{O}_{Y}(1-l)) = 0$ for $l \geq 2$. For $i = 1$, dualizing the Euler exact sequence, twisting by $\mathcal{O}_Y(1-l)$ and taking cohomology, we have that  
$$0 \to H^{0}(\mathcal{Q}_Y^*(1-l)) \to H^{0}(\mathcal{O}_Y(1-l))^{\oplus 7} $$
Since $H^{0}(\mathcal{O}_Y(1-l)) = 0$ for $l \geq 2$, we have that  $H^0(\mathcal{Q}_Y^*(1-l)) = 0$ for $l \geq 2$. 

\color{black}
We now show that $H^1(\mathcal{Q}_Y^*(1-l)) = 0$ for $l \geq 3$. Dualizing the exact sequence

$$0 \to \mathcal{U} \to \mathcal{O}_{G_5}^{\oplus 7} \to \mathcal{Q} \to 0$$

restricting to $Y$ and tensoring with $\mathcal{O}_Y(1-l)$ we have 

$$0 \to H^0(\mathcal{U}_Y^*(1-l)) \to H^1(\mathcal{Q}_Y^*(1-l)) \to 0$$

since $H^0(\mathcal{O}_Y(1-l)) = H^1(\mathcal{O}_Y(1-l)) = 0$. So it is enough to prove that $H^0(\mathcal{U}_Y^*(1-l)) = 0$. As before, we consider the exact sequence 
$$0 \to \mathcal{U}^* \otimes \mathcal{I}_Y (1-l) \to \mathcal{U}^*(1-l) \to \mathcal{U}_Y^*(1-l) \to 0$$

Taking cohomology, we see that it is enough to show the vanishing of $H^0(\mathcal{U}^*(1-l))$ and $H^1(\mathcal{U}^* \otimes \mathcal{I}_Y(1-l))$. \par

We have $H^0(\mathcal{U}^*(1-l)) = H^0(\Sigma^{(2-l,1-l)} \mathcal{U}^*) = 0 $ for $l \geq 3$.
To show the second vanishing, we tensor the Koszul resolution of the ideal sheaf $\mathcal{I}_Y$ by $\mathcal{U}^*(1-l)$ to have \\

\begin{multline*}
0 \to \bigwedge^7 E^* \otimes \mathcal{U}^*(1-l) \to \bigwedge^6 E^* \otimes \mathcal{U}^*(1-l) \to \bigwedge^5 E^* \otimes \mathcal{U}^*(1-l) \to \bigwedge^4 E^* \otimes \mathcal{U}^*(1-l) \to \bigwedge^3 E^* \otimes \mathcal{U}^*(1-l) \to \bigwedge^2 E^* \otimes \mathcal{U}^*(1-l) \\ \to E^* \otimes \mathcal{U}^*(1-l)  \to \mathcal{U}^* \otimes \mathcal{I}_Y (1-l)) \to 0
\end{multline*}

We need to show that $H^{k}(\bigwedge^k E^* \otimes \mathcal{U}^*(1-l)) = 0$ for $ 1 \leq k \leq 7$. We have $$\bigwedge^k E^* = \bigwedge^{k}(\mathcal{Q})(-k) \oplus (\bigwedge^{k-1}(\mathcal{Q})(-(k-1)) \otimes \mathcal{O}_{G_5}(-1))^{\oplus 2} \oplus (\bigwedge^{k-2}(\mathcal{Q})(-(k-2)) \otimes \mathcal{O}_{G_5}(-2)) $$

We look at the terms $H^{k}(\mathcal{U}^*(1-l) \otimes \bigwedge^{k}(\mathcal{Q})(-k))$, $H^{k}(\mathcal{U}^*(1-l) \otimes \bigwedge^{k-1}(\mathcal{Q})(-(k-1)) \otimes \mathcal{O}_{G_5}(-1))$ and $H^{k}(\mathcal{U}^*(1-l) \otimes \bigwedge^{k-2}(\mathcal{Q})(-(k-2)) \otimes \mathcal{O}_{G_5}(-2))$. We first look at $H^{k}(\mathcal{U}^*(1-l) \otimes \bigwedge^{k}(\mathcal{Q})(-k))$. \\
We have,
\begin{multline*}
H^{k}(\mathcal{U}^*(1-l) \otimes \bigwedge^{k}(\mathcal{Q})(-k)) = H^{k}(\mathcal{U}^*(-k-l+1) \otimes \bigwedge^k \mathcal{Q}) \\ = H^{k}(\Sigma^{(-k-l+2,-k-l+1)} \mathcal{U}^* \otimes \Sigma^{\underbrace{(1,..,1}_{\text{$k$}},0,0,..)}\mathcal{Q}) \\ = H^{k}(\Sigma^{(-k-l+2,-k-l+1)} \mathcal{U}^* \otimes \Sigma^{(0,0,..0,\underbrace{-1,..,-1}_{\text{$k$}})}\mathcal{Q}^*)  \\
= H^{k}(L_{(-k-l+2, -k-l+1, 0,..0, \underbrace{-1,..,-1}_{\text{$k$}}))}) 
\end{multline*}
 
One can now check as before that this is zero for $1 \leq k \leq 7$, and $l \geq 3$. In fact it is enough to check this term for $1 \leq k \leq 5$. 

Now we look at $H^{k}(\mathcal{U}^*(1-l) \otimes \bigwedge^{k-1}(\mathcal{Q})(-(k-1)) \otimes \mathcal{O}_{G_5}(-1))$. Once again we have,  that 
\begin{multline*}
   H^{k}(\mathcal{U}^*(1-l) \otimes \bigwedge^{k-1}(\mathcal{Q})(-(k-1)) \otimes \mathcal{O}_{G_5}(-1)) = H^{k}(\mathcal{U}^*(-k-l+1) \otimes \bigwedge^{k-1} \mathcal{Q})) \\ 
   = H^{k}(\Sigma^{(-k-l+2,-k-l+1)} \mathcal{U}^* \otimes \Sigma^{\underbrace{(1,..,1}_{\text{$k-1$}},0,0,..)}\mathcal{Q}) \\ = H^{k}(\Sigma^{(-k-l+2,-k-l+1)} \mathcal{U}^* \otimes \Sigma^{(0,0,..0,\underbrace{-1,..,-1}_{\text{$k-1$}})}\mathcal{Q}^*)
  \\ = H^{k}(L_{(-k-l+2, -k-l+1, 0,..0, \underbrace{-1,..,-1}_{\text{$k-1$}})}) 
\end{multline*}

One can now check \textcolor{black}{that the above is} zero for $1 \leq k \leq 7$, and $l \geq 3$. In fact it is enough to check this term for $1 \leq k \leq 6$.

We look at $H^{k}(\mathcal{U}^*(1-l) \otimes \bigwedge^{k-2}(\mathcal{Q})(-(k-2)) \otimes \mathcal{O}_{G_5}(-2))$. Once again we have,  that 
\begin{multline*}
   H^{k}(\mathcal{U}^*(1-l) \otimes \bigwedge^{k-2}(\mathcal{Q})(-(k-2)) \otimes \mathcal{O}_{G_5}(-2)) = H^{k}(\mathcal{U}^*(-k-l+1) \otimes \bigwedge^{k-2} \mathcal{Q})) \\ 
   = H^{k}(\Sigma^{(-k-l+2,-k-l+1)} \mathcal{U}^* \otimes \Sigma^{\underbrace{(1,..,1}_{\text{$k-2$}},0,0,..)}\mathcal{Q}) \\
  = H^{k}(\Sigma^{(-k-l+2,-k-l+1)} \mathcal{U}^* \otimes \Sigma^{0,0,..,\underbrace{-1,..,-1}_{\text{$k-2$}})}\mathcal{Q}^*) \\
   = H^{k}(L_{(-k-l+2, -k-l+1, 0,..0, \underbrace{-1,..,-1}_{\text{$k-2$}})}) 
\end{multline*}

One can now check as before that both factors are zero for $1 \leq k \leq 7$, and $l \geq 2$. This concludes proof for $H^1(\mathcal{Q}_Y^*(1-l)) = 0$.

\color{black}

Hence we have
$$0 \to H^i(T_Y(-l)) \to H^i(T_{G_5}|_Y(-l))$$ for $l \geq 2$. As before, we have
$$0 \to T_{G_5} \otimes \mathcal{I}_Y (-l) \to T_{G_5}(-l) \to T_{G_5}|_Y(-l) \to 0$$
Taking cohomology, we have
$$H^i(T_{G_5}(-l)) \to H^i(T_{G_5}|_Y(-l)) \to H^{i+1}(T_{G_5} \otimes \mathcal{I}_Y (-l))$$

We have that $H^i(T_{G_5}(-l)) = H^{10-i}(\Omega_{G_5}(l-7)) = 0$. To show that $H^{i+1}(T_{G_3} \otimes \mathcal{I}_Y (-l)) = 0$, using the Koszul resolution of the ideal sheaf $\mathcal{I}_Y$, we have \\
\begin{multline*}
0 \to \bigwedge^7 E^* \otimes T_{G_5}(-l) \to \bigwedge^6 E^* \otimes T_{G_5}(-l) \to \bigwedge^5 E^* \otimes T_{G_5}(-l) \to \bigwedge^4 E^* \otimes T_{G_5}(-l) \to \bigwedge^3 E^* \otimes T_{G_5}(-l) \to \bigwedge^2 E^* \otimes T_{G_5}(-l) \\ \to E^* \otimes T_{G_5}(-l)  \to T_{G_5} \otimes \mathcal{I}_Y (-l)) \to 0
\end{multline*}

We need to show that $H^{k+i}(\bigwedge^k E^* \otimes T_{G_5}(-l)) = 0$ for $ 1 \leq k \leq 7$. We have $$\bigwedge^k E^* = \bigwedge^{k}(\mathcal{Q})(-k) \oplus (\bigwedge^{k-1}(\mathcal{Q})(-(k-1)) \otimes \mathcal{O}_{G_5}(-1))^{\oplus 2} \oplus (\bigwedge^{k-2}(\mathcal{Q})(-(k-2)) \otimes \mathcal{O}_{G_5}(-2)) $$

We look at the terms $H^{k+i}(T_{G_5}(-l) \otimes \bigwedge^{k}(\mathcal{Q})(-k))$, $H^{k+i}(T_{G_5}(-l) \otimes \bigwedge^{k-1}(\mathcal{Q})(-(k-1)) \otimes \mathcal{O}_{G_5}(-1))$ and $H^{k+i}(T_{G_5}(-l) \otimes \bigwedge^{k-2}(\mathcal{Q})(-(k-2)) \otimes \mathcal{O}_{G_5}(-2))$. We first look at $H^{k+i}(T_{G_5}(-l) \otimes \bigwedge^{k}(\mathcal{Q})(-k))$. \\
We have by Pieri's rule,
\begin{multline*}
H^{k+i}(T_{G_5}(-l) \otimes \bigwedge^{k}(\mathcal{Q})(-k)) = H^{k+i}(\mathcal{U}^*(-k-l) \otimes \mathcal{Q} \otimes \bigwedge^k \mathcal{Q}) \\ = H^{k+i}(\Sigma^{(-k-l+1,-k-l)} \mathcal{U}^* \otimes \Sigma^{\underbrace{(1,..,1}_{\text{$k+1$}},0,0,..)}\mathcal{Q}) \oplus H^{k+i}(\Sigma^{(-k-l+1,-k-l)} \mathcal{U}^* \otimes \Sigma^{\underbrace{(2,..,1}_{\text{$k$}},0,0,..)}\mathcal{Q}) \\
= H^{k+i}(\Sigma^{(-k-l+1,-k-l)} \mathcal{U}^* \otimes \Sigma^{(0,0,..,\underbrace{-1,..,-1}_{\text{$k+1$}})}\mathcal{Q}^*) \oplus H^{k+i}(\Sigma^{(-k-l+1,-k-l)} \mathcal{U}^* \otimes \Sigma^{(0,0,..,\underbrace{-2,..,-1}_{\text{$k$}})}\mathcal{Q}^*) \\= H^{k+i}(L_{(-k-l+1, -k-l, 0,0,.., \underbrace{-1,..,-1}_{\text{$k+1$}})}) \oplus 
H^{k+i}(L_{(-k-l+1, -k-l, 0,0,.., \underbrace{-2,..,-1}_{\text{$k$}})})
\end{multline*}
 
 One can now check as before that both factors are zero for $1 \leq k \leq 7$, $i = 1,2$ and $l \geq 2$. In fact it is enough to check this term for $1 \leq k \leq 5$.

Now we look at $H^{k+i}(T_{G_5}(-l) \otimes \bigwedge^{k-1}(\mathcal{Q})(-(k-1)) \otimes \mathcal{O}_{G_5}(-1))$. Once again we have by Pieri's rule,  that 
\begin{multline*}
   H^{k+i}(T_{G_5}(-l) \otimes \bigwedge^{k-1}(\mathcal{Q})(-(k-1)) \otimes \mathcal{O}_{G_5}(-1)) = H^{k+i}(\mathcal{U}^*(-k-l) \otimes \mathcal{Q} \otimes \bigwedge^{k-1} \mathcal{Q})) \\ 
   = H^{k+i}(\Sigma^{(-k-l+1,-k-l)} \mathcal{U}^* \otimes \Sigma^{\underbrace{(1,..,1}_{\text{$k$}},0,0,..)}\mathcal{Q}) \oplus H^{k+i}(\Sigma^{(-k-l+1,-k-l)} \mathcal{U}^* \otimes \Sigma^{\underbrace{(2,..,1}_{\text{$k-1$}},0,0,..)}\mathcal{Q}) \\ = H^{k+i}(\Sigma^{(-k-l+1,-k-l)} \mathcal{U}^* \otimes \Sigma^{(0,0,..,\underbrace{-1,..,-1}_{\text{$k$}})}\mathcal{Q}^*) \oplus H^{k+i}(\Sigma^{(-k-l+1,-k-l)} \mathcal{U}^* \otimes \Sigma^{0,0,..,\underbrace{-2,..,-1}_{\text{$k-1$}})}\mathcal{Q}^*) \\
   = H^{k+i}(L_{(-k-l+1, -k-l, 0,0,.., \underbrace{-1,..,-1}_{\text{$k$}})}) \oplus 
H^{k+i}(L_{(-k-l+1, -k-l, 0,0,.., \underbrace{-2,..,-1}_{\text{$k-1$}})})
\end{multline*}

 One can now check as before that both factors are zero for $1 \leq k \leq 7$, $i = 1,2$ and $l \geq 2$. In fact it is enough to check this term for $1 \leq k \leq 6$.

We look at $H^{k+i}(T_{G_5}(-l) \otimes \bigwedge^{k-2}(\mathcal{Q})(-(k-2)) \otimes \mathcal{O}_{G_5}(-2))$. Once again we have by Pieri's rule,  that 
\begin{multline*}
   H^{k+i}(T_{G_5}(-l) \otimes \bigwedge^{k-2}(\mathcal{Q})(-(k-2)) \otimes \mathcal{O}_{G_5}(-2)) = H^{k+i}(\mathcal{U}^*(-k-l) \otimes \mathcal{Q} \otimes \bigwedge^{k-2} \mathcal{Q})) \\ 
   = H^{k+i}(\Sigma^{(-k-l+1,-k-l)} \mathcal{U}^* \otimes \Sigma^{\underbrace{(1,..,1}_{\text{$k-1$}},0,0,..)}\mathcal{Q}) \oplus H^{k+i}(\Sigma^{(-k-l+1,-k-l)} \mathcal{U}^* \otimes \Sigma^{\underbrace{(2,..,1}_{\text{$k-2$}},0,0,..)}\mathcal{Q}) \\ = H^{k+i}(\Sigma^{(-k-l+1,-k-l)} \mathcal{U}^* \otimes \Sigma^{(0,0,..,\underbrace{-1,..,-1}_{\text{$k-1$}})}\mathcal{Q}^*) \oplus H^{k+i}(\Sigma^{(-k-l+1,-k-l)} \mathcal{U}^* \otimes \Sigma^{(0,0,..,\underbrace{-2,..,-1}_{\text{$k-2$}})}\mathcal{Q}^*) \\
   = H^{k+i}(L_{(-k-l+1, -k-l,0,0,.., \underbrace{-1,..,-1}_{\text{$k-1$}})}) \oplus 
H^{k+i}(L_{(-k-l+1, -k-l, 0,0,.., \underbrace{-2,..,-1}_{\text{$k-2$}})})
\end{multline*}

One can now check as before that both factors are zero for $1 \leq k \leq 7$, $i = 1,2$ and $l \geq 2$. This concludes our proof for family $1.9$. 

\textcolor{black}{\bf{The proof for family $1.8$}}. In this case, we have that $Y$ is the zero locus of a section of the vector bundle $E = \bigwedge^2 \mathcal{U}^* \oplus \mathcal{O}_G(1)^{\oplus 3}$ on the Grassmannian $G = \textrm{Gr}(3,6)$. We have $K_G = \mathcal{O}_G(-6)$. Note that $\mathcal{U}^*$ has rank $3$ and hence $\bigwedge^2 \mathcal{U}^* = \mathcal{U} \otimes \textrm{det}(\mathcal{U}^*) = \mathcal{U}(1)$. We have that $N_{Y/G} = E|_Y = \mathcal{U}_Y(1) \oplus \mathcal{O}_Y(1)^{\oplus 3}$ and hence $\textrm{det}(N_{Y/G}) = \mathcal{O}_Y(5)$. Hence $K_Y = \mathcal{O}_{Y}(-1)$. 
Twisting the sequence 
$$0 \to T_Y \to T_{G}|_Y \to \mathcal{U}_Y(1) \oplus \mathcal{O}_Y(1)^{\oplus 3} \to 0$$ 
Twisting by $\mathcal{O}_Y(-l)$ and taking cohomology, we have for $i = 1,2$
$$H^{i-1}(\mathcal{U}_Y(1-l) \oplus H^{i-1}(\mathcal{O}_{Y}(1-l))^{\oplus 3}) \to H^i(T_Y(-l)) \to H^i(T_{G}|_Y(-l))$$
Clearly, $H^0(\mathcal{O}_{Y}(1-l)) = 0$ for $l \geq 2$. Dualizing the Euler exact sequence, twisting by $\mathcal{O}_Y(1-l)$ and taking cohomology, we have that  
$$0 \to H^0(\mathcal{U}_Y(1-l)) \to H^0(\mathcal{O}_Y(1-l))^{\oplus 6} $$
Considering that $H^0(\mathcal{O}_Y(1-l)) = 0$, we have that $H^0(\mathcal{U}_Y(1-l)) = 0$. \\
\color{black}
Let us now show $H^1(\mathcal{U}_Y(1-l)) = 0$. Note that using the Euler exact sequence, we have 

$$H^0(\mathcal{O}_Y(1-l))^{\oplus 6} \to H^0(\mathcal{Q}_Y(1-l)) \to H^1(\mathcal{U}_Y(1-l)) \to H^1(\mathcal{O}_Y(1-l))^{\oplus 6}$$

Considering that both the flanking terms are zero we are left to show the vanishing of $H^0(\mathcal{Q}_Y(1-l))$. Now as before using the defining exact sequence of $Y$ inside the Grassmannian, we it is enough to show the vanishing of $H^0(\mathcal{Q}(1-l))$ and $H^1(\mathcal{Q} \otimes \mathcal{I}_Y(1-l))$. The former can be checked by Borel-Weil-Bott vanishing theorem for $l \geq 3$. For the latter we use the resolution of the ideal sheaf and tensor by $\mathcal{Q}(1-l)$ to get

\begin{multline*}
 0 \to \bigwedge^6 E^* \otimes \mathcal{Q}(1-l) \to \bigwedge^5 E^* \otimes \mathcal{Q}(1-l) \to \bigwedge^4 E^* \otimes \mathcal{Q}(1-l) \to \bigwedge^3 E^* \otimes \mathcal{Q}(1-l) \to \bigwedge^2 E^* \otimes \mathcal{Q}(1-l) \\ \to E^* \otimes \mathcal{Q}(1-l)  \to \mathcal{Q} \otimes \mathcal{I}_Y (1-l)) \to 0
\end{multline*}

It is enough to show that $H^k(\mathcal{Q}(1-l) \otimes \bigwedge^k E^*) = 0$ for $1 \leq k \leq 6$. We have $$\bigwedge^k E^* = \bigwedge^{k}(\mathcal{U}^*)(-k) \oplus (\bigwedge^{k-1}(\mathcal{U}^*)(-(k-1)) \otimes \mathcal{O}_{G}(-1))^{\oplus 3} \oplus (\bigwedge^{k-2}(\mathcal{U}^*)(-(k-2)) \otimes \mathcal{O}_{G}(-2))^{\oplus 3} \oplus (\bigwedge^{i-3}(\mathcal{U}^*)(-(k-3)) \otimes \mathcal{O}_{G}(-3)) $$

We first look at $H^{k}(\mathcal{Q}(1-l) \otimes \bigwedge^{k}(\mathcal{U}^*)(-k))$. We have 
\begin{multline*}
    H^{k}(\mathcal{Q}(1-l) \otimes \bigwedge^{k}(\mathcal{U}^*)(-k)) = H^{k}(\bigwedge^k (\mathcal{U}^*)(-k-l+1) \otimes \mathcal{Q}) \\
    = H^{k}(\Sigma^{\underbrace{(-k-l+2,..,-k-l+2}_{\text{k}}, \underbrace{-k-l+1,..,-k-l+1}_{\text{3-k}})} \mathcal{U}^* \otimes \Sigma^{(1,0,0)} \mathcal{Q}) \\
    = H^{k}(\Sigma^{\underbrace{(-k-l+2,..,-k-l+2}_{\text{k}}, \underbrace{-k-l+1,..,-k-l+1}_{\text{3-k}})} \mathcal{U}^* \otimes \Sigma^{(0,0,-1)} \mathcal{Q}^*)
\end{multline*}

One can now check that this is zero for $1 \leq k \leq 6$ for $l \geq 3$ (Actually it is enough to check this for $1 \leq k \leq 3$).

Now we look at $H^{k}(\mathcal{Q}(1-l) \otimes \bigwedge^{k-1}(\mathcal{U}^*)(-k))$. We have 
\begin{multline*}
    H^{k}(\mathcal{Q}(1-l) \otimes \bigwedge^{k-1}(\mathcal{U}^*)(-k)) = H^{k}(\bigwedge^{k-1} (\mathcal{U}^*)(-k-l+1) \otimes \mathcal{Q}) \\
    = H^{k}(\Sigma^{\underbrace{(-k-l+2,..,-k-l+2}_{\text{k-1}}, \underbrace{-k-l+1,..,-k-l+1}_{\text{4-k}})} \mathcal{U}^* \otimes \Sigma^{(1,0,0)} \mathcal{Q})
    \\ = H^{k}(\Sigma^{\underbrace{(-k-l+2,..,-k-l+2}_{\text{k-1}}, \underbrace{-k-l+1,..,-k-l+1}_{\text{4-k}})} \mathcal{U}^* \otimes \Sigma^{(0,0,-1)} \mathcal{Q}^*)
\end{multline*}

One can now check that this is zero for $1 \leq k \leq 6$ for $l \geq 3$ (It is enough to check this for $1 \leq k \leq 4$).

Now we look at $H^{k}(\mathcal{Q}(1-l) \otimes \bigwedge^{k-1}(\mathcal{U}^*)(-k))$. We have 
\begin{multline*}
    H^{k}(\mathcal{Q}(1-l) \otimes \bigwedge^{k-2}(\mathcal{U}^*)(-k)) = H^{k}(\bigwedge^{k-2} (\mathcal{U}^*)(-k-l+1) \otimes \mathcal{Q}) \\
    = H^{k}(\Sigma^{\underbrace{(-k-l+2,..,-k-l+2}_{\text{k-2}}, \underbrace{-k-l+1,..,-k-l+1}_{\text{5-k}})} \mathcal{U}^* \otimes \Sigma^{(1,0,0)} \mathcal{Q})
    \\ = H^{k}(\Sigma^{\underbrace{(-k-l+2,..,-k-l+2}_{\text{k-2}}, \underbrace{-k-l+1,..,-k-l+1}_{\text{5-k}})} \mathcal{U}^* \otimes \Sigma^{(0,0,-1)} \mathcal{Q}^*)
\end{multline*}

One can now check that this is zero for $1 \leq k \leq 6$ for $l \geq 3$ (It is enough to check this for $2 \leq k \leq 5$).

Finally we look at $H^{k}(\mathcal{Q}(1-l) \otimes \bigwedge^{k-3}(\mathcal{U}^*)(-k))$. We have 
\begin{multline*}
    H^{k}(\mathcal{Q}(1-l) \otimes \bigwedge^{k-2}(\mathcal{U}^*)(-k)) = H^{k}(\bigwedge^{k-3} (\mathcal{U}^*)(-k-l+1) \otimes \mathcal{Q}) \\
    = H^{k}(\Sigma^{\underbrace{(-k-l+2,..,-k-l+2}_{\text{k-3}}, \underbrace{-k-l+1,..,-k-l+1}_{\text{6-k}})} \mathcal{U}^* \otimes \Sigma^{(1,0,0)} \mathcal{Q})
    \\ = H^{k}(\Sigma^{\underbrace{(-k-l+2,..,-k-l+2}_{\text{k-3}}, \underbrace{-k-l+1,..,-k-l+1}_{\text{6-k}})} \mathcal{U}^* \otimes \Sigma^{(0,0,-1)} \mathcal{Q}^*)
\end{multline*}

One can now check that this is zero for $1 \leq k \leq 6$ for $l \geq 3$ (It is enough to check this for $3 \leq k \leq 6$).

\color{black}

Hence we have for $i = 1,2$,
$$0 \to H^i(T_Y(-l)) \to H^i(T_{G}|_Y(-l))$$ for $l \geq 2$ and $i = 1,2$. 

We have
$$0 \to T_{G} \otimes \mathcal{I}_Y (-l) \to T_{G}(-l) \to T_{G}|_Y(-l) \to 0$$
Taking cohomology, we have
$$H^i(T_{G}(-l)) \to H^i(T_{G}|_Y(-l)) \to H^{i+1}(T_{G} \otimes \mathcal{I}_Y (-l))$$
We have that $H^i(T_{G}(-l)) = H^{9-i}(\Omega_{G}(l-6)) = 0$. To show that $H^{i+1}(T_{G} \otimes \mathcal{I}_Y (-l)) = 0$, using the Koszul resolution of the ideal sheaf $\mathcal{I}_Y$ we have
 \begin{multline*}
 0 \to \bigwedge^6 E^* \otimes T_{G}(-l) \to \bigwedge^5 E^* \otimes T_{G}(-l) \to \bigwedge^4 E^* \otimes T_{G}(-l) \to \bigwedge^3 E^* \otimes T_{G}(-l) \to \bigwedge^2 E^* \otimes T_{G}(-l) \\ \to E^* \otimes T_{G}(-l)  \to T_{G} \otimes \mathcal{I}_Y (-l)) \to 0
\end{multline*}

We need to show that $H^{k+i}(\bigwedge^k E^* \otimes T_{G}(-l)) = 0$ for $1 \leq k \leq 6$. We have $$\bigwedge^k E^* = \bigwedge^{k}(\mathcal{U}^*)(-k) \oplus (\bigwedge^{k-1}(\mathcal{U}^*)(-(k-1)) \otimes \mathcal{O}_{G}(-1))^{\oplus 3} \oplus (\bigwedge^{k-2}(\mathcal{U}^*)(-(k-2)) \otimes \mathcal{O}_{G}(-2))^{\oplus 3} \oplus (\bigwedge^{i-3}(\mathcal{U}^*)(-(k-3)) \otimes \mathcal{O}_{G}(-3)) $$
We first look at $H^{k+i}(T_G(-l) \otimes \bigwedge^{k}(\mathcal{U}^*)(-k))$. We have by Pieri's rule
\begin{multline*}
    H^{k+i}(T_G(-l) \otimes \bigwedge^{k}(\mathcal{U}^*)(-k)) = H^{k+i}(\mathcal{U}^* \otimes \bigwedge^k \mathcal{U}^*(-k-l) \otimes \mathcal{Q}) \\
    = H^{k+i}(\Sigma^{\underbrace{(-k-l+1,..,-k-l+1}_{\text{k+1}}, \underbrace{-k-l,..,-k-l}_{\text{2-k}})} \mathcal{U}^* \otimes \Sigma^{(1,0,0)} \mathcal{Q}) \\ \oplus H^{k+i}(\Sigma^{(-k-l+2,\underbrace{-k-l+1,..,-k-l+1}_{\text{k-1}}, \underbrace{-k-l,..,-k-l}_{\text{3-k}})} \mathcal{U}^* \otimes \Sigma^{(1,0,0)} \mathcal{Q}) \\
    = H^{k+i}(\Sigma^{\underbrace{(-k-l+1,..,-k-l+1}_{\text{k+1}}, \underbrace{-k-l,..,-k-l}_{\text{2-k}})} \mathcal{U}^* \otimes \Sigma^{(0,0,-1)} \mathcal{Q}^*) \\ \oplus H^{k+i}(\Sigma^{(-k-l+2,\underbrace{-k-l+1,..,-k-l+1}_{\text{k-1}}, \underbrace{-k-l,..,-k-l}_{\text{3-k}})} \mathcal{U}^* \otimes \Sigma^{(0,0,-1)} \mathcal{Q}^*)
\end{multline*}

One can check that both the terms are equal to zero for $1 \leq k \leq 6$ (Actually it is enough to check the first term for $1 \leq k \leq 2$ and the second term for $1 \leq k \leq 3$) for $l \geq 2$.

Now we look at $H^{k+i}(T_G(-l) \otimes \bigwedge^{k-1}(\mathcal{U}^*)(-k))$. We have by Pieri's rule
\begin{multline*}
    H^{k+i}(T_G(-l) \otimes \bigwedge^{k-1}(\mathcal{U}^*)(-k)) = H^{k+i}(\mathcal{U}^* \otimes \bigwedge^{k-1} \mathcal{U}^*(-k-l) \otimes \mathcal{Q}) \\
    = H^{k+i}(\Sigma^{\underbrace{(-k-l+1,..,-k-l+1}_{\text{k}}, \underbrace{-k-l,..,-k-l}_{\text{3-k}})} \mathcal{U}^* \otimes \Sigma^{(1,0,0)} \mathcal{Q}) \\ \oplus H^{k+i}(\Sigma^{(-k-l+2,\underbrace{-k-l+1,..,-k-l+1}_{\text{k-2}}, \underbrace{-k-l,..,-k-l}_{\text{4-k}})} \mathcal{U}^* \otimes \Sigma^{(1,0,0)} \mathcal{Q})
    \\ = H^{k+i}(\Sigma^{\underbrace{(-k-l+1,..,-k-l+1}_{\text{k}}, \underbrace{-k-l,..,-k-l}_{\text{3-k}})} \mathcal{U}^* \otimes \Sigma^{(0,0,-1)} \mathcal{Q}^*) \\ \oplus H^{k+i}(\Sigma^{(-k-l+2,\underbrace{-k-l+1,..,-k-l+1}_{\text{k-2}}, \underbrace{-k-l,..,-k-l}_{\text{4-k}})} \mathcal{U}^* \otimes \Sigma^{(0,0,-1)} \mathcal{Q}^*)
\end{multline*}

One can check that both the terms are equal to zero for $1 \leq k \leq 6$ (It is enough to check the first term for $1 \leq k \leq 3$ and the second term for $2 \leq k \leq 4$) for $l \geq 2$.

Next we look at $H^{k+i}(T_G(-l) \otimes \bigwedge^{k-2}(\mathcal{U}^*)(-k))$. We have by Pieri's rule
\begin{multline*}
    H^{k+i}(T_G(-l) \otimes \bigwedge^{k-2}(\mathcal{U}^*)(-k)) = H^{k+i}(\mathcal{U}^* \otimes \bigwedge^{k-2} \mathcal{U}^*(-k-l) \otimes \mathcal{Q}) \\
    = H^{k+i}(\Sigma^{\underbrace{(-k-l+1,..,-k-l+1}_{\text{k-1}}, \underbrace{-k-l,..,-k-l}_{\text{4-k}})} \mathcal{U}^* \otimes \Sigma^{(1,0,0)} \mathcal{Q}) \\ \oplus H^{k+i}(\Sigma^{(-k-l+2,\underbrace{-k-l+1,..,-k-l+1}_{\text{k-3}}, \underbrace{-k-l,..,-k-l}_{\text{5-k}})} \mathcal{U}^* \otimes \Sigma^{(1,0,0)} \mathcal{Q}) \\
    =  H^{k+i}(\Sigma^{\underbrace{(-k-l+1,..,-k-l+1}_{\text{k-1}}, \underbrace{-k-l,..,-k-l}_{\text{4-k}})} \mathcal{U}^* \otimes \Sigma^{(0,0,-1)} \mathcal{Q}^*) \\ \oplus H^{k+i}(\Sigma^{(-k-l+2,\underbrace{-k-l+1,..,-k-l+1}_{\text{k-3}}, \underbrace{-k-l,..,-k-l}_{\text{5-k}})} \mathcal{U}^* \otimes \Sigma^{(0,0,-1)} \mathcal{Q}^*)
\end{multline*}

One can check that both the terms are equal to zero for $1 \leq k \leq 6$ (It is enough to check the first term for $1 \leq k \leq 4$ and the second term for $3 \leq k \leq 5$) for $l \geq 2$.

\vspace{0.5cm}

Finally, we look at $H^{k+i}(T_G(-l) \otimes \bigwedge^{k-3}(\mathcal{U}^*)(-k))$. We have by Pieri's rule
\begin{multline*}
    H^{k+i}(T_G(-l) \otimes \bigwedge^{k-3}(\mathcal{U}^*)(-k)) = H^{k+i}(\mathcal{U}^* \otimes \bigwedge^{k-3} \mathcal{U}^*(-k-l) \otimes \mathcal{Q}) \\
    = H^{k+i}(\Sigma^{\underbrace{(-k-l+1,..,-k-l+1}_{\text{k-2}}, \underbrace{-k-l,..,-k-l}_{\text{5-k}})} \mathcal{U}^* \otimes \Sigma^{(1,0,0)} \mathcal{Q}) \\ \oplus H^{k+i}(\Sigma^{(-k-l+2,\underbrace{-k-l+1,..,-k-l+1}_{\text{k-4}}, \underbrace{-k-l,..,-k-l}_{\text{6-k}})} \mathcal{U}^* \otimes \Sigma^{(1,0,0)} \mathcal{Q}) \\ 
    = H^{k+i}(\Sigma^{\underbrace{(-k-l+1,..,-k-l+1}_{\text{k-2}}, \underbrace{-k-l,..,-k-l}_{\text{5-k}})} \mathcal{U}^* \otimes \Sigma^{(0,0,-1)} \mathcal{Q}^*) \\ \oplus H^{k+i}(\Sigma^{(-k-l+2,\underbrace{-k-l+1,..,-k-l+1}_{\text{k-4}}, \underbrace{-k-l,..,-k-l}_{\text{6-k}})} \mathcal{U}^* \otimes \Sigma^{(0,0,-1)} \mathcal{Q}^*)
\end{multline*}

One can check that both the terms are equal to zero for $1 \leq k \leq 5$ (It is enough to check the first term for $2 \leq k \leq 5$ and the second term for $4 \leq k \leq 6$) for $l \geq 2$.

\vspace{0.5cm}

\textcolor{black}{\bf{ The proof for family $1.10$}}. In this case, we have that $Y$ is the zero locus of a section of the vector bundle $E = (\bigwedge^2 \mathcal{U}^*)^{\oplus 3}$ on the Grassmannian $G = \textrm{Gr}(3,7)$. We have $K_G = \mathcal{O}_G(-7)$. Note that $\mathcal{U}^*$ has rank $3$ and hence $\bigwedge^2 \mathcal{U}^* = \mathcal{U} \otimes \textrm{det}(\mathcal{U}^*) = \mathcal{U}(1)$. We have that $N_{Y/G} = E|_Y = \mathcal{U}_Y(1)^{\oplus 3}$ and hence $\textrm{det}(N_{Y/G}) = \mathcal{O}_Y(6)$. Hence $K_Y = \mathcal{O}_{Y}(-1)$. 
Twisting the sequence 
$$0 \to T_Y \to T_{G}|_Y \to \mathcal{U}_Y(1)^ {\oplus 3} \to 0$$ 
Twisting by $\mathcal{O}_Y(-l)$ and taking cohomology, we have 
$$H^{i-1}(\mathcal{U}_Y(1-l)^{\oplus 3} \to H^i(T_Y(-l)) \to H^i(T_{G}|_Y(-l))$$
Dualizing the Euler exact sequence, twisting by $\mathcal{O}_Y(1-l)$ and taking cohomology, we have that  
$$0 \to H^0(\mathcal{U}_Y(1-l)) \to H^0(\mathcal{O}_Y(1-l))^{\oplus 7} $$
Considering that $H^0(\mathcal{O}_Y(1-l)) = 0$ for $l \geq 2$, we have that $H^0(\mathcal{U}_Y(1-l)) = 0$. The vanishing of $H^0(\mathcal{U}_Y(1-l)) = 0$ follows as in the proof of family $1.8$. 
Hence we have
$$0 \to H^1(T_Y(-l)) \to H^1(T_{G}|_Y(-l))$$ for $l \geq 2$. As before, we have
$$0 \to T_{G} \otimes \mathcal{I}_Y (-l) \to T_{G}(-l) \to T_{G}|_Y(-l) \to 0$$
Taking cohomology, we have
$$H^1(T_{G}(-l)) \to H^1(T_{G}|_Y(-l)) \to H^2(T_{G} \otimes \mathcal{I}_Y (-l))$$
We have that $H^1(T_{G_5}(-l)) = H^{11}(\Omega_{G_5}(l-7)) = 0$. To show that $H^2(T_{G_3} \otimes \mathcal{I}_Y (-l)) = 0$, using the Koszul resolution of the ideal sheaf $\mathcal{I}_Y$, we have \\

$$0 \to \bigwedge^9 E^* \otimes T_{G_5}(-l) \to \cdot \cdot \cdot \to E^* \otimes T_{G_5}(-l)  \to T_{G_5} \otimes \mathcal{I}_Y (-l)) \to 0$$

We need to show that $H^{i+1}(\bigwedge^i E^* \otimes T_{G_5}(-l)) = 0$ for $ 1 \leq i \leq 9$. Note that $E^* = \mathcal{U}^*(-1)^{\oplus 3}$. Hence 
$$\bigwedge^i E^* = \bigoplus_{\substack{m+n+k = i \\ 0 \leq m,n,k \leq 3}}\bigwedge^m (\mathcal{U}^*(-1)) \otimes \bigwedge^n (\mathcal{U}^*(-1)) \otimes \bigwedge^k (\mathcal{U}^*(-1)) = \bigoplus_{\substack{m+n+k = i \\ 0 \leq m,n,k \leq 3}}(\bigwedge^m \mathcal{U}^* \otimes \bigwedge^n \mathcal{U}^* \otimes \bigwedge^k \mathcal{U}^*)(-i)$$
So we need to check the vanishing of $ H^{i+1}(\mathcal{U}^*(-l-i) \otimes \bigwedge^m \mathcal{U}^* \otimes \bigwedge^n \mathcal{U}^* \otimes \bigwedge^k \mathcal{U}^* \otimes \mathcal{Q}) $ for $1 \leq i \leq 9$, $m+n+k = i$ and $0 \leq m,n,k \leq 3$. We have to show the vanishing for the tuples $(i, m, n, k) = (9,3,3,3), (8,3,3,2), (7,3,3,1), (7,3,2,2), (6,3,3,0), (6,3,2,1), (6,2,2,2), (5,3,2,0), (5,3,1,1), (5,2,2,1), (4,3,1,0), (4,2,2,0),$ $ (4,2,1,1), (3,3,0,0), (3,2,1,0), (3,1,1,1), (2,2,0,0), (2,1,1,0), (1,1,0,0) $.
We work out the case $i=6, m=2, n=2, k=2$. The rest of the cases are easier and follow either as this case or as one of the cases worked out before. We have by Pieri's rule

\begin{align*}
 H^7(\mathcal{U}^*(-l-6) \otimes (\bigwedge^2\mathcal{U}^*)^{\otimes 3} \otimes \mathcal{Q}) & = H^7(\mathcal{U}^*(-l-3) \otimes \mathcal{U}^{\otimes 3} \otimes \mathcal{Q}) \\
& = H^7(\mathcal{U}^*(-l-3) \otimes \Sigma^{(2,1,0)}\mathcal{U} \otimes \mathcal{Q})^{\oplus 2}  \bigoplus H^7(\mathcal{U}^*(-l-3) \otimes \Sigma^{(3,0,0)}\mathcal{U} \otimes \mathcal{Q}) \\
& \bigoplus H^7(\mathcal{U}^*(-l-3) \otimes \Sigma^{(1,1,1)}\mathcal{U} \otimes \mathcal{Q}) \\
& = H^7(\mathcal{U}^*(-l-3) \otimes \Sigma^{(0,-1,-2)}\mathcal{U}^* \otimes \mathcal{Q})^{\oplus 2}  \bigoplus H^7(\mathcal{U}^*(-l-3) \otimes \Sigma^{(0,0, -3)}\mathcal{U}^* \otimes \mathcal{Q}) \\
& \bigoplus H^7(\mathcal{U}^*(-l-3) \otimes \Sigma^{(-1,-1,-1)}\mathcal{U}^* \otimes \mathcal{Q}) \\
& = H^7(\Sigma^{(1,-1,-2)}\mathcal{U}^*(-l-3) \otimes \mathcal{Q})^{\oplus 2} \bigoplus H^7(\Sigma^{(0,0,-2)}\mathcal{U}^*(-l-3) \otimes \mathcal{Q})^{\oplus 2} \\
&  \bigoplus H^7(\Sigma^{(0,-1,-1)}\mathcal{U}^*(-l-3) \otimes \mathcal{Q})^{\oplus 2} \bigoplus H^7(\Sigma^{(1,0,-3)}\mathcal{U}^*(-l-3) \otimes \mathcal{Q}) \\
&  \bigoplus H^7(\Sigma^{(0,0,-2)}\mathcal{U}^*(-l-3) \otimes \mathcal{Q}) \bigoplus H^7(\Sigma^{(0,-1,-1)}\mathcal{U}^*(-l-3) \otimes \mathcal{Q}) \\
& = H^7(\Sigma^{(-l-2,-l-4,-l-5)}\mathcal{U}^* \otimes \mathcal{Q})^{\oplus 2} \bigoplus H^7(\Sigma^{(-l-3,-l-3,-l-5)}\mathcal{U}^* \otimes \mathcal{Q})^{\oplus 2} \\ 
& \bigoplus H^7(\Sigma^{(-l-3,-l-4,-l-4)}\mathcal{U}^* \otimes \mathcal{Q}) 
 \bigoplus H^7(\Sigma^{(-l-2,-l-3,-l-6)}\mathcal{U}^* \otimes \mathcal{Q}) \\
& \bigoplus H^7(\Sigma^{(-l-3,-l-3,-l-5)}\mathcal{U}^* \otimes \mathcal{Q}) \bigoplus H^7(\Sigma^{(-l-3,-l-4,-l-4)}\mathcal{U}^* \otimes \mathcal{Q}) \\
& = H^7(\Sigma^{(-l-2,-l-4,-l-5)}\mathcal{U}^* \otimes \Sigma^{(0,0,-1)}\mathcal{Q}^*)^{\oplus 2} \bigoplus H^7(\Sigma^{(-l-3,-l-3,-l-5)}\mathcal{U}^* \otimes \Sigma^{(0,0,-1)}\mathcal{Q}^*)^{\oplus 2} \\ 
& \bigoplus H^7(\Sigma^{(-l-3,-l-4,-l-4)}\mathcal{U}^* \otimes \Sigma^{(0,0,-1)}\mathcal{Q}^*) 
 \bigoplus H^7(\Sigma^{(-l-2,-l-3,-l-6)}\mathcal{U}^* \otimes \Sigma^{(0,0,-1)}\mathcal{Q}^*) \\
& \bigoplus H^7(\Sigma^{(-l-3,-l-3,-l-5)}\mathcal{U}^* \otimes \Sigma^{(0,0,-1)}\mathcal{Q}^*) \bigoplus H^7(\Sigma^{(-l-3,-l-4,-l-4)}\mathcal{U}^* \otimes \Sigma^{(0,0,-1)}\mathcal{Q}^*)
\end{align*}
Each of the last six terms can be checked to be equal to zero for $l \geq 2$.

\vspace{0.5cm}

\textcolor{black}{\textbf{The proof for family $1.1$ and $1.12$ are very much same. We only prove for family $1.12$.}} $Y$ is a double cover of $\mathbb{P}^3$ branched along a smooth quartic $B$. Let $\pi: Y \to \mathbb{P}^3$ denote the double cover. Then $K_Y = \pi^*(\mathcal{O}_{\mathbb{P}^3}(-2))$. Considering the exact sequence 
$$0 \to T_Y \to \pi^*T_{\mathbb{P}^3} \to N_{\pi} \to 0$$

Tensoring with $\pi^*(\mathcal{O}_{\mathbb{P}^3}(-l))$ and taking cohomology we have for $i = 1,2$

$$H^{i-1}(N_{\pi} \otimes \pi^*(\mathcal{O}_{\mathbb{P}^3}(-l)) \to H^i(T_Y \otimes \pi^*(\mathcal{O}_{\mathbb{P}^3}(-l))) \to H^i(\pi^*T_{\mathbb{P}^3} \otimes \pi^*(\mathcal{O}_{\mathbb{P}^3}(-l)))$$

Note that by pushing forward, we have that 
$H^{i-1}(N_{\pi} \otimes \pi^*(\mathcal{O}_{\mathbb{P}^3}(-l)) = H^{i-1}(\mathcal{O}_{\mathbb{P}^3}(4-l)|_B)$ \textcolor{black}{(see \cite[Lemma $2.5$]{GGP10})}. For $i = 1$, this is equal to $h^0(\mathcal{O}_{\mathbb{P}^3}(1)) = 4$ for $l = 3$, $1$ for $l = 4$ and $0$ for $l \geq 5$. For $i = 2$, this is equal to $0$. Now note that
$H^i(\pi^*T_{\mathbb{P}^3} \otimes \pi^*(\mathcal{O}_{\mathbb{P}^3}(-l)) = H^i(T_{\mathbb{P}^3} \otimes (\mathcal{O}_{\mathbb{P}^3}(-l))) \oplus H^i(T_{\mathbb{P}^3} \otimes (\mathcal{O}_{\mathbb{P}^3}(-l-2)))$

Tensoring the Euler sequence 

$$0 \to \mathcal{O}_{\mathbb{P}^3} \to \mathcal{O}_{\mathbb{P}^3}(1)^{\oplus 4} \to T_{\mathbb{P}^3} \to 0$$

by $\mathcal{O}_{\mathbb{P}^3}(-k)$ and taking cohomology we have that 

$$H^i(\mathcal{O}_{\mathbb{P}^3}(1-k)) \to H^i(T_{\mathbb{P}^3} \otimes (\mathcal{O}_{\mathbb{P}^3}(-k))) \to H^{i+1}(\mathcal{O}_{\mathbb{P}^3}(-k)) \to H^{i+1}(\mathcal{O}_{\mathbb{P}^3}(1-k))^4 $$

For $k \geq 2$, $i = 1$, $H^1(T_{\mathbb{P}^3} \otimes (\mathcal{O}_{\mathbb{P}^3}(-k))) = 0$. For $i = 2$, note that the right flanking map is 
$$ H^{3}(\mathcal{O}_{\mathbb{P}^3}(-k)) \to H^{3}(\mathcal{O}_{\mathbb{P}^3}(1-k)) \otimes H^0({\mathcal{O}_{\mathbb{P}^3}(1))} $$
which by Serre-duality is the dual of the map 
$$H^{0}(\mathcal{O}_{\mathbb{P}^3}(k-5)) \otimes H^0({\mathcal{O}_{\mathbb{P}^3}(1))} \to H^{0}(\mathcal{O}_{\mathbb{P}^3}(k-4))$$
which is surjective for $k \geq 5$ and hence the dual map is injective. So for $k \geq 5$, $H^2(T_{\mathbb{P}^3} \otimes (\mathcal{O}_{\mathbb{P}^3}(-k))) = 0$.  
\vspace{0.5cm}

\textcolor{black}{\bf{The proof for families $1.13$, $1.14$, $1.16$ and $1.17$ are very similar. We only prove for family $1.13$.}} In this case, $Y$ is a hypersurface of degree $3$ in $\mathbb{P}^4$. $K_Y = \mathcal{O}_Y(-2)$ where $\mathcal{O}_Y(1)$ is the class of a hyperplane section. Then tensoring the sequence

$$0 \to T_Y \to T_{\mathbb{P}^4}|_Y \to \mathcal{O}_Y(3) \to 0$$

by $\mathcal{O}_Y(-l)$ and taking cohomology, we have that for $ i = 1,2$,

$$H^{i-1}(\mathcal{O}_Y(3-l)) \to H^i(T_Y(-l)) \to H^i(T_{\mathbb{P}^4}|_Y(-l)) $$

We have $h^0(\mathcal{O}_Y(3-l)) = 1$ for $l = 3$ and $0$ if $l \geq 4$. Further $h^1(\mathcal{O}_Y(3-l)) = 0$ for any $l$. Now tensoring the pullback of the Euler sequence to $Y$ by $\mathcal{O}_Y(-l)$, we have

$$H^i(\mathcal{O}_Y(1-l))^{\oplus 5} \to H^i(T_{\mathbb{P}^4}|_Y(-l)) \to H^{i+1}(\mathcal{O}_Y(-l))$$
Since for $l \geq 3$ and $i = 1$, both flanking terms are zero, we have that $h^1(T_Y(-l)) \leq 1$ for $l = 3$ and $0$ for $l \geq 4$. On the other hand for $ i = 2$, the last long exact sequence is 
$$0 \to H^2(T_{\mathbb{P}^4}|_Y(-l)) \to H^{3}(\mathcal{O}_Y(-l)) \to H^3(\mathcal{O}_Y(1-l)) \otimes H^0(\mathcal{O}_{\mathbb{P}^4}(1))$$
The right flanking map is the dual of the map 

$$ H^0(\mathcal{O}_Y(l-3)) \otimes H^0(\mathcal{O}_{\mathbb{P}^4}(1)) \to H^0(\mathcal{O}_Y(l-2)) $$

which is surjective for $l \geq 3$. Hence the right flanking map is injective. for $l \geq 3$. Hence $H^2(T_Y(-l)) = 0$ for $l \geq 3$.

\vspace{0.5cm}

\textcolor{black}{\bf{The proof for families $1.2$, $1.3$ and $1.4$ are very similar.}} In each of these cases $Y$ is a complete intersection in $\mathbb{P}^N$ for some $N$ and $K_Y = \mathcal{O}_Y(-1)$. We have the sequence

$$0 \to T_Y \to T_{\mathbb{P}^N}|_Y \to N_{Y/\mathbb{P}^N} \to 0$$

Tensoring with $\mathcal{O}_Y(-l)$ and taking cohomology, we have that for $ i = 1,2$, 

$$H^{i-1}(N_{Y/\mathbb{P}^N}(-l)) \to H^i(T_Y(-l)) \to H^i(T_{\mathbb{P}^N}|_Y(-l)) $$

In each case one can see that $H^1(N_{Y/\mathbb{P}^N}(-l)) = 0$ and using the pullback of the Euler sequence to $Y$ one can see that $H^1(T_{\mathbb{P}^N}|_Y(-l)) = 0$ for $l \geq 2$. In case $1.2$, $H^0(N_{Y/\mathbb{P}^N}(-l)) = H^0(\mathcal{O}_Y(4-l))$ and hence we have that $H^0(\mathcal{O}_Y(4-l)) = 1$ if $l = 4$ and $0$ if $l \geq 5$. In case $1.3$, $H^0(N_{Y/\mathbb{P}^N}(-l)) = H^0(\mathcal{O}_Y(2-l)) \oplus H^0(\mathcal{O}_Y(3-l))$ which is $1$ if $l = 3$ and $0$ if $l \leq 4$. In case $1.4$, $H^0(N_{Y/\mathbb{P}^N}(-l)) = H^0(\mathcal{O}_Y(2-l))^{\oplus 3} $ which is $3$ if $l = 2$ and $0$ if $l \leq 3$. Finally the vanishing of $H^2(T_{\mathbb{P}^N}|_Y(-l))$ for $l \geq 2$ follows from the projective normality of complete intersections as before.

\vspace{0.5cm}

\textcolor{black}{\bf{Let us prove for family $1.11$}}. In this case, $Y$ is a degree $6$ hypersurface inside $\mathbb{P}(1,1,1,2,3) = \mathbb{P}$. $K_Y = \mathcal{O}_Y(-2)$. We have as usual, after tensoring by $\mathcal{O}_Y(-l)$, the sequence

$$0 \to T_Y(-l) \to T_{\mathbb{P}}|_Y(-l) \to \mathcal{O}_Y(6-l) \to 0 $$

We have that for $l=4$, $h^0(\mathcal{O}_Y(6-l)) = h^0(\mathcal{O}_Y(2)) = 11$, for $l = 5$, $h^0(\mathcal{O}_Y(6-l)) = h^0(\mathcal{O}_Y(1)) = 3$, for $l = 6$, $h^0(\mathcal{O}_Y(6-l)) = h^0(\mathcal{O}_Y) = 1$,  while if $l \geq 7$, $h^0(\mathcal{O}_Y(6-l)) = 0$. Further for any value of $l$, $h^1(\mathcal{O}_Y(6-l)) = 0$.  Now note that if $\overline{\Omega}_{\mathbb{P}}$ denotes the sheaf of regular differential forms as defined in \cite[Section $(2)$]{Dol82}, then $T_{\mathbb{P}} = \overline{\Omega}_{\mathbb{P}}^*$. To see this, further note that  if $i: \mathbb{P}_{\textrm{sm}} \hookrightarrow \mathbb{P}$ denotes the inclusion of smooth locus, then $i_*(T_{\mathbb{P}_{\textrm{sm}}}) = T_{\mathbb{P}}$ since both are reflexive sheaves that agree on the smooth locus which has codimension at least $2$. Further it is shown in \cite[$2.2.4$]{Dol82}, that $i_*(\Omega_{\mathbb{P}_{\textrm{sm}}}) = \overline{\Omega_{\mathbb{P}}}$. We have $i_*(\Omega_{\mathbb{P}_{\textrm{sm}}}^*) = i_*(\Omega_{\mathbb{P}_{\textrm{sm}}})^*$ because both are reflexive sheaves that agree on the smooth locus which has codimension at least $2$. So finally
$T_{\mathbb{P}} = i_*(T_{\mathbb{P}_{\textrm{sm}}}) = i_*(\Omega_{\mathbb{P}_{\textrm{sm}}}^*) = i_*(\Omega_{\mathbb{P}_{\textrm{sm}}})^*  = \overline{\Omega}_{\mathbb{P}}^*$. Now we dualize the Euler sequence 

$$0 \to \overline{\Omega}_{\mathbb{P}} \to \mathcal{O}_{\mathbb{P}}(-1)^{\oplus 3} \oplus \mathcal{O}_{\mathbb{P}}(-2) \oplus \mathcal{O}_{\mathbb{P}}(-3) \to \mathcal{O}_{\mathbb{P}} \to 0$$

to get 

$$0 \to \mathcal{O}_{\mathbb{P}} \to  \mathcal{O}_{\mathbb{P}}(1)^{\oplus 3} \oplus \mathcal{O}_{\mathbb{P}}(2) \oplus \mathcal{O}_{\mathbb{P}}(3) \to T_{\mathbb{P}} \to  0$$

Now restricting to $Y$, tensoring by $\mathcal{O}_Y(-l)$ and taking cohomology we have by \cite[$1.4$]{Dol82}, that $H^1(T_{\mathbb{P}}(-l)|_Y) = 0$ for any $l$. \textcolor{black}{Now we show that $H^2(T_{\mathbb{P}}(-l)|_Y) = 0$ for $l \geq 6$. To see this first note that by restricting the Euler exact sequence to $Y$, twisting by $\mathcal{O}_Y(-l)$ and taking the cohomology, we have the exact sequence
$$0 \to H^2(T_{\mathbb{P}}|_Y(-l)) \to H^3(\mathcal{O}_Y(1-l))^{\oplus 3} \oplus H^3(\mathcal{O}_Y(2-l)) \oplus H^3(\mathcal{O}_Y(3-l)) \xrightarrow{\lambda} H^3(\mathcal{O}_Y(-l))$$
So the result follows if we show that the map $\lambda$ is injective or equivalently the dual $\lambda^*$ is surjective. Note that by Serre-duality $\lambda^*$ is the map between, 
$$H^0(\mathcal{O}_Y(l-3))^{\oplus 3} \oplus H^0(\mathcal{O}_Y(l-4)) \oplus H^0(\mathcal{O}_Y(l-5)) \xrightarrow{\lambda^*} H^0(\mathcal{O}_Y(l-2))$$
We show that the map $H^0(\mathcal{O}_Y(l-3))^{\oplus 3} \to H^0(\mathcal{O}_Y(l-2))$ is already surjective. Note that $H^0(\mathcal{O}_Y(l-3))^{\oplus 3} = H^0(\mathcal{O}_Y(l-3)) \otimes H^0(\mathcal{O}_{\mathbb{P}}(1))$ (from the construction of the Euler sequence). Now the fact that the above map is surjective follows from the fact that the map $H^0(\mathcal{O}_{\mathbb{P}}(l-3)) \otimes H^0(\mathcal{O}_{\mathbb{P}}(1)) \to H^0(\mathcal{O}_{\mathbb{P}}(l-2))$ is surjective( see \cite{NO02}) and the fact that the maps of global sections $H^0(\mathcal{O}_{\mathbb{P}}(l-3)) \to H^0(\mathcal{O}_{Y}(l-3))$ and $H^0(\mathcal{O}_{\mathbb{P}}(l-2)) \to H^0(\mathcal{O}_{Y}(l-2))$ are both surjective.  }
\QEDA

\subsection{Invariants of the Hilbert components of precisely $1-$ extendable canonical surfaces}\label{invariants of hyperplane sections}
The relative position of the invariants invariants $(\chi, K^2)$ of the Hilbert components of the precisely $1-$ extendable canonical surfaces that we obtain with respect to the fundamental inequalities in the geography of surfaces of general type is shown as follows.

\begin{figure}[h!]
    \centering
    \includegraphics[scale = 0.1]{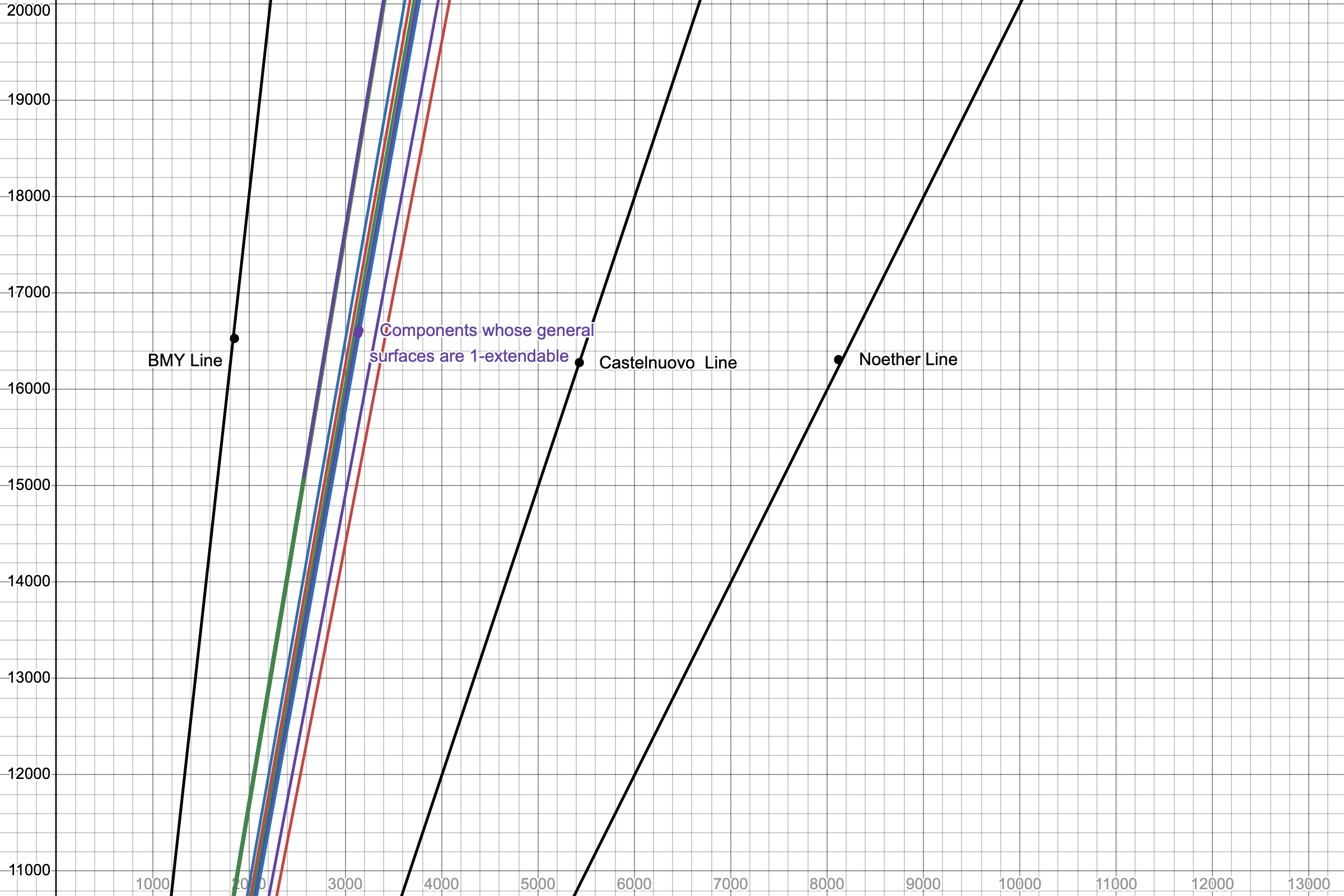}
    \label{geography}
\end{figure}

More specifically, the invariants are as follows :

\begin{itemize}
    \item[(1)] For families $1.1-1.10$, i.e, when $Y$ is of index $1$, a surface in $\mathcal{S}^Y_l$ has $$(\chi,K^2) = (\displaystyle\frac{1}{3}(-K_Y^3)l^3+(\displaystyle\frac{1}{6}(-K_Y^3)+4)l, 2l^3(-K_Y^3))$$ 
    \begin{figure}[h!]
    \centering
    \includegraphics[scale = 0.1]{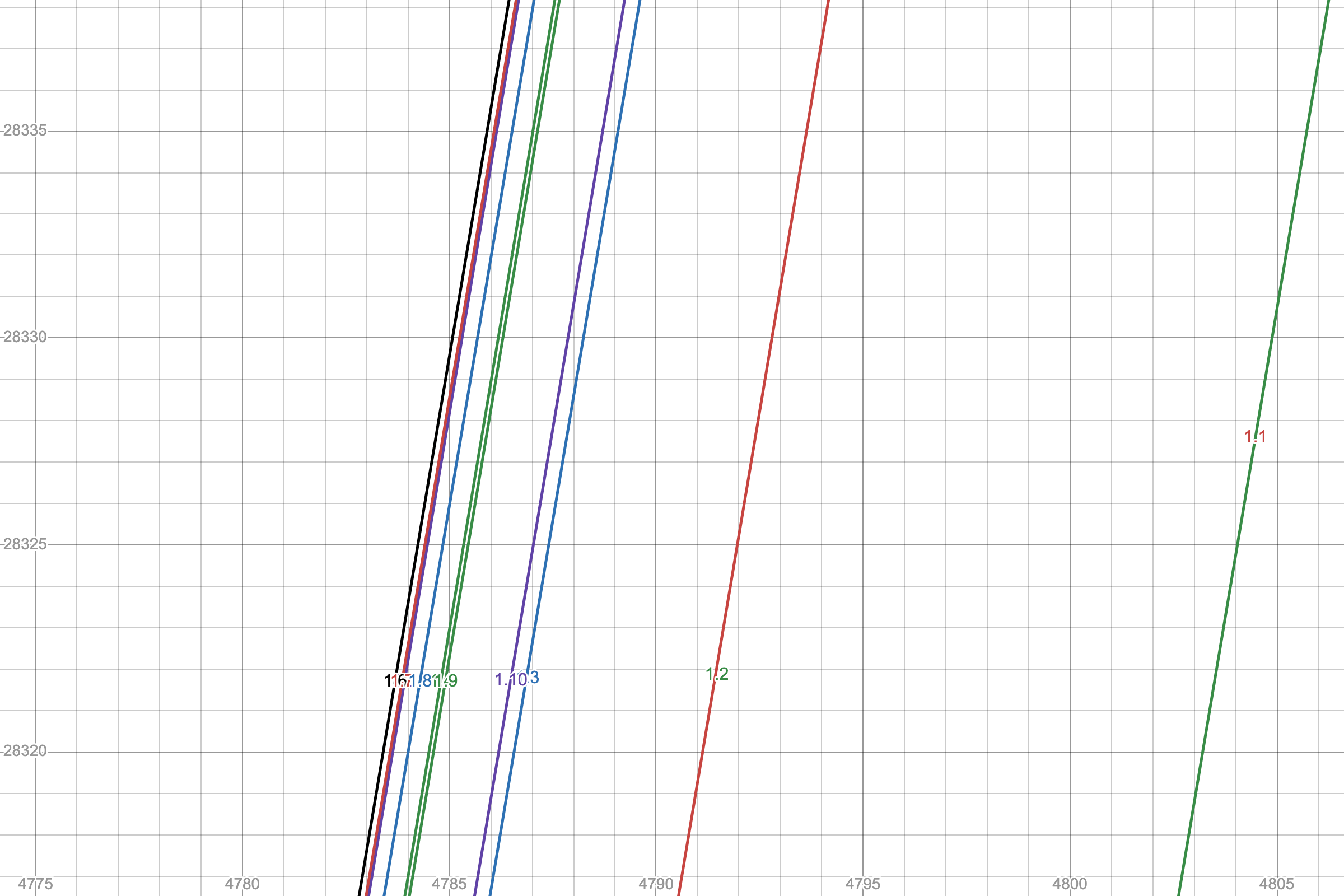}
    \label{index 1}
    \end{figure}
   
   
   \item[(2)] For families $1.11-1.17$, i.e, when $Y$ is of index $j$, a surface in $\mathcal{S}^Y_l$ has $$(\chi,K^2) = (\displaystyle\frac{1}{3}H^3l^3+ \displaystyle\frac{j-1}{2}H^3l^2+(\displaystyle\frac{1}{6}H^3+2)l+1, 2H^3l^3)$$ where $H$ is the generator of the Picard group of $Y$.
    \begin{figure}[h!]
    \centering
    \includegraphics[scale = 0.1]{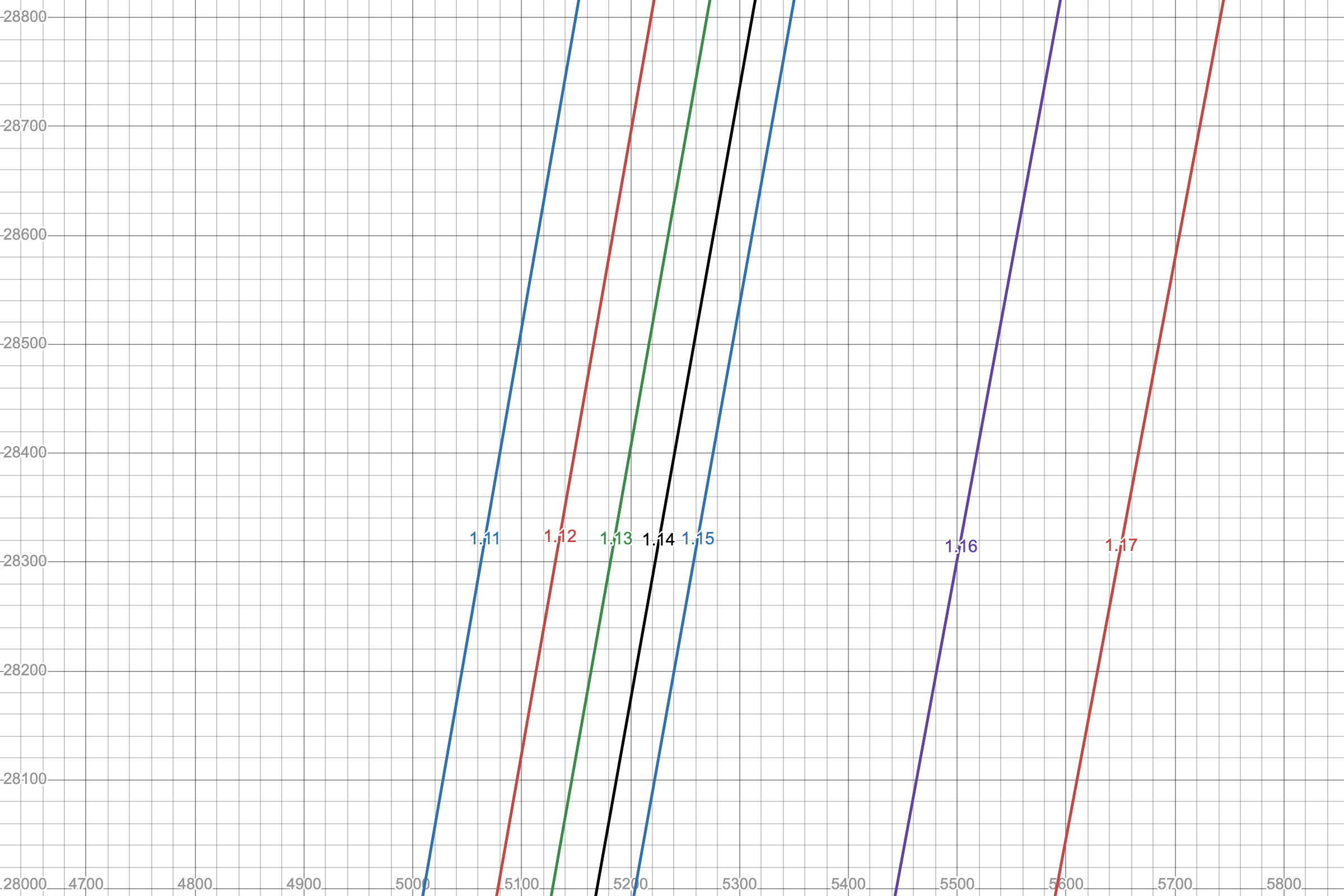}
    \label{index > 1}
    \end{figure}

\end{itemize}


\bibliographystyle{plain}

\end{document}